\documentclass[12pt]{amsart}

\usepackage{amssymb, amsmath, amsthm}
\usepackage{xcolor}
\definecolor{myred}{rgb}{0.2,0,0}
\definecolor{myblue}{rgb}{0,0,0.6}
\definecolor{mygreen}{rgb}{0,0.2,0}
\usepackage[colorlinks=true,linkcolor=myred,citecolor=myblue,urlcolor=mygreen]{hyperref}
\usepackage[centering]{geometry}

\newcommand\multsum{\mathop{\sum\cdots\sum}\limits}

\newcommand{\les}{\leqslant}
\newcommand{\ges}{\geqslant}

\newtheorem{theorem}{Theorem}
\newtheorem{lemma}{Lemma}

\newtheorem{corollary}{Corollary}

\newtheorem{remark}{Remark}

\begin{document}
	
\title{Fractional parts of non-integer powers of primes. II} 
\author{Andrei Shubin}
\address{Department of Mathematics, Caltech, 1200 E. California Blvd., Pasadena, CA, 91125, USA}
\email{\href{mailto:ashubin@caltech.edu}{ashubin@caltech.edu}}
		
\maketitle

\begin{abstract}
We continue to study the distribution of prime numbers $p$, satisfying the condition $\{ p^{\alpha} \} \in I \subset [0; 1)$, in arithmetic progressions. In the paper, we prove an analogue of Bombieri-Vinogradov theorem for $0 < \alpha < 1/9$ with the level of distribution $\theta = 2/5 - (3/5) \alpha$, which improves the previous result corresponding to $\theta \les 1/3$.
\end{abstract}

\section{Introduction}
\label{sec1}
	
As in the previous work~\cite{23} let us denote by $\mathbb{E} \subset \mathbb{N}$ the subsequence of natural numbers
$$
	\bigl\{n \in \mathbb{N} : \{ n^{\alpha}  \} \in I \bigr\},
$$ where $\alpha > 0$ is any fixed non-integer, $I$ is any subinterval of $[0; 1)$. The distribution of primes from $\mathbb{E}$ was studied by a number of authors including Vinogradov, Linnik, Kaufman, Gritsenko, Balog, Harman, Tolev and many others (see \cite{1, 2, 4, 5, 9, 10, 11, 12, 17, 18, 19, 22, 25, 26, 27}). One of the main results of this investigation is the asymptotic formula for the proportion of such primes:
$$
	\sum_{\substack{p \les X \\ p \in \mathbb{E}}} 1 = |I| \cdot \pi(X) + O\bigl(X^{1-\vartheta(\alpha)}\bigr),
$$ where the exponent $0 < \vartheta(\alpha) < 1$ had been sharpening for different values of $\alpha$ until at least 2006. For $0 < \alpha < 1$ see \cite{2, 10, 12, 19, 22, 26}. The case $\alpha > 1$, $\alpha \notin \mathbb{N}$ is covered in \cite{1, 4, 5, 9, 18, 27}. 

The other direction concerns the existence of infinite number of primes from a very thin subset of integers of the form $\bigl\{ n \in \mathbb{N}: \{ \sqrt{n} \} < n^{-c + \varepsilon}  \bigr\}$ for fixed $c > 0$ and arbitrary small $\varepsilon > 0$. The first such result is due to Vinogradov, who proved this for all $c \les 1/10$ (see~\cite[Ch.~4]{28}). Later it was improved in the work of Kaufman~\cite{17} to all $c \les 0.16310\ldots$ and in the unpublished work of Harman to $c \les 0.2139\ldots$. In~\cite{17} Kaufman also showed that the Riemann Hypothesis implies this result for all $c \les 1/4$. It was proved for all $c \les 1/4$ unconditionally in the papers of Balog~\cite{2} and Harman~\cite{12}. Finally, Harman and Lewis established the result for all $c \les 0.262$ in~\cite{13}.
	
We focus on the distribution of $p \in \mathbb{E}$ in the arithmetic progressions of the form $qn + a$, $(a,q)=1$. An analogue of Bombieri-Vinogradov theorem for such subset is usually given by the	 inequality
\begin{equation} \label{BV}
	\sum_{q \les Q} \max_{(a,q)=1} \biggl| \sum_{\substack{p \les X \\ p \in \mathbb{E} \\ p \equiv a\pmod{q}}} 1 - \frac{1}{\varphi(q)} \sum_{\substack{p \les X \\ p \in \mathbb{E}}} 1 \biggr| \ll_A \frac{X}{(\log X)^A},
\end{equation} where $A > 0$ can be arbitrarily large, $\varepsilon > 0$ is arbitrarily small, $Q = X^{\theta - \varepsilon}$ and $\theta$ is called the ``level of distribution''. 
	
Tolev~\cite{25} showed~\eqref{BV} for $\alpha = 1/2$ and all $\theta \les 1/4$. Later Gritsenko and Zinchenko~\cite{11} extended this result to all $1/2 \les \alpha < 1$ and $\theta \les 1/3$. In~\cite{23} the author further extended it to all $\alpha > 0$, $\alpha \notin \mathbb{N}$ and $\theta \les 1/3$. In the present paper we improve this result for small $\alpha$, namely we show~\eqref{BV} for all $\theta \les 2/5 - (3/5) \alpha$, which goes beyond the range $\theta \les 1/3$ if $0 < \alpha < 1/9$. 

The proof is based on the estimation of the exponential sum of the form
$$
	\sum_{\substack{X \les p < 2X \\ p \equiv a\pmod{q}}} e\bigl( hp^{\alpha} \bigr),
$$ where $e(x) := e^{2\pi ix}$. The desired upper bound is given by

\begin{theorem} \label{thm1}
	Suppose that $0 < \alpha < 1/9$ is fixed non-integer, $\theta, \varepsilon, C$ are fixed constants satisfying the conditions $0 < \varepsilon < \alpha / 100$, $\varepsilon < \theta < 2/5 - (3/5)\alpha$, $C \ges 1$, and suppose that $1 \les  h \les (\log X)^C$, $2 < q \les X^{\theta - \varepsilon}$, $1 \les a \les q-1$, $(a,q)=1$. Then the sum
	$$
		T = \sum_{\substack{X \les  p < 2X  \\ p \equiv a\pmod{q}}} e\bigl( h p^{\alpha}\bigr) 
	$$ satisfies the estimate
	\begin{equation} \label{to_prove}
		T \ll \frac{X}{q} (\log X)^{-A}
	\end{equation} with an arbitrarily large $A > 0$.
\end{theorem}

\begin{corollary} \label{cor1}
	Let $0 < \alpha < 1/9$ be a fixed non-integer, $\varepsilon > 0$ is arbitrary small number. Then for any $q \les X^{2/5 - (3/5)\alpha - \varepsilon}$, $a$, $(a,q)=1$, and any given subinterval $I \subset [0; 1)$ the following asymptotic formula holds true:
	$$
		\pi_{I} (X; q, a) := \sum_{\substack{p \les X \\ \{ p^{\alpha} \} \in I \\ p \equiv a\pmod{q}}} 1 = |I| \cdot \pi(X; q,a) + O \biggl( \frac{\pi(X; q,a)}{(\log X)^A} \biggr)
	$$ with any fixed $A > 0$.
\end{corollary}

The next corollary is the analogue of Bombieri-Vinogradov theorem:

\begin{corollary} \label{cor2}
	Let $0 < \alpha < 1/9$ be fixed, $I = [c; d) \subset [0; 1)$, $\mathbb{E} = \bigl\{ n \in \mathbb{N} : \{ n^{\alpha} \} \in I \bigr\}$, and let $\theta, \varepsilon, A$ be fixed constants such that $0 < \varepsilon < \theta < 2/5 - (3/5)\alpha$, $\varepsilon < \alpha / 100$, $A > 0$. Next, let $2 < Q \les X^{\theta - \varepsilon}$. Then the following inequality holds true:
	$$
		\sum_{q \les Q} \max_{(a,q)=1} \biggl| \sum_{\substack{p \les X \\ p \in \mathbb{E} \\ p \equiv a\pmod{q}}} 1 - \frac{1}{\varphi(q)} \sum_{\substack{p \les X \\ p \in \mathbb{E}}} 1 \biggr| \ll \frac{X}{(\log X)^A}.
	$$
\end{corollary} 

\begin{remark}
	In the analogue of Bombieri-Vinogradov theorem one can apparently go beyond the level of distirubtion $\theta = 2/5 - (3/5)\alpha$ using the large sieve. This is a work in progress.
\end{remark}

The way one deduces corollaries \hyperref[cor1]{1} and \hyperref[cor2]{2} from \hyperref[thm1]{Theorem~1} is explicated in the Section 2 of~\cite{23}. In this work we only focus on the proof of \hyperref[thm1]{Theorem~1}. The main difference between this proof and the proof of Theorem 1 from~\cite{23} is the application of Heath-Brown identity~\cite{14} in place of Vaughan identity~\cite[Ch.~13]{15}. A key new ingredient is a combinatorial decomposition of the initial sum
$$
	W := \sum_{\substack{X \les n < Y \\ n \equiv a\pmod{q}}} \Lambda(n) e\bigl( hn^{\alpha} \bigr), \qquad X < Y \les 2X,
$$ into sums of three types. This decomposition is given by \hyperref[lemma1]{Lemma~1} in \hyperref[sec2]{Section~2} (see also \cite[Lemma~3.1]{21}). Then $W$ splits into
\begin{gather*}
	W_I = \sum_{m \sim M} a_m \sum_{\substack{n \sim N \\ mn \equiv a (\text{mod} \ q)}} f(n) e \bigl(h (mn)^{\alpha}\bigr), \\
	W_{II} = \sum_{m \sim M} \beta_m \sum_{\substack{n \sim N \\ mn \equiv a (\text{mod} \ q)}} \gamma_n e \bigl(h (mn)^{\alpha}\bigr), \\
	W_{III} = \sum_{m \sim M} f_1 (m) \sum_{n \sim N} f_2 (n) \sum_{\substack{k \sim K \\ mnk \equiv a \pmod{q}}} f_3(k) e\bigl( h (mnk)^{\alpha} \bigr), \\
\end{gather*} where $a_m, \beta_m$, $\gamma_n$ are real coefficients, $f, f_1, f_2, f_3$ are smooth functions and ``$x \sim X$'' means $X\Theta^{-1} \les x \les X\Theta$ for some fixed $\Theta$, $1 < \Theta < 2$. All of the coefficients $|a_m|, |\beta_m|, |\gamma_n|, |f|, |f_1|, |f_2|, |f_3|$ do not exceed $X^{\varepsilon}$. Since the number of sums of each type is bounded by a constant, we have
$$
	W \ll \bigl|W_I\bigr| + \bigl|W_{II}\bigr| + \bigl|W_{III}\bigr|.
$$ We treat $W_I, W_{II}, W_{III}$ separately. 

The estimation of type I and type II sums is very similar to the one in~\cite{23}, and it only requires the classical van der Corput second derivative test~\cite[Ch.~1, Theorem~5]{16} due to the small size of $\alpha$. This estimation is carried out in \hyperref[sec3]{Section~3} and \hyperref[sec4]{Section~4}.

The estimation of type III sum is contained in \hyperref[sec5]{Section~5}. It is a little more delicate. In this case all three variables are roughly of the same size $K \approx N \approx M \approx X^{1/3}$. This implies that for large values of $q$, precisely, $q \ges X^{1/3}$, the sum over $k, mnk \equiv a\pmod{q}$ might be empty or contain only one term, so one cannot get a cancellation from the inner sum. This is the main reason why the previous method does not work in this case. The new idea is to remove the congruence condition $mnk \equiv a\pmod{q}$ using orthogonality of Dirichlet characters and then apply Poisson summation formula to any two of three sums over $m,n,k$. In this way we replace them by two shorter sums, and the new expression for $W_{III}$ would have the form
$$
	W_{III} = \sum_{m \sim M} f_4 (m) \sum_{u \sim U} \sum_{v \sim V} f_5 (uv) S_q (u,v),
$$ where $UV \ll KN$, $f_4 (m)$ and $f_5 (uv)$ will be specified in the end of \hyperref[sec5]{Section~5}. We then estimate $|f_4 (m)|$ and $|f_5 (uv)|$ trivially and apply Weil's bound $|S_q (u,v)| \les \sqrt{q} \tau(q) (u,v,q)^{1/2}$ for Kloosterman sum (see, for example,~\cite{6}) to obtain the upper bound for $|W_{III}|$.

The application of Poisson summation requires a smoothed sum. So before estimating $|W_I|$, $|W_{II}|$ and $|W_{III}|$ we would slightly adjust $W$ in \hyperref[sec2]{Section~2}. Precisely, we remove the sharp bounds for $m,n,k$ using the smooth partition of unity, which is also described in~\cite[Section~3]{21}. To deal with the oscillating integrals arising after Poisson summation we apply the method of stationery phase. The necessary tools are given by \hyperref[lemma2]{Lemma~2} and \hyperref[lemma3]{Lemma~3} in \hyperref[sec5]{Section~5}. They are proven in~\cite{3}.

\section{Initial steps. Heath-Brown identity. Smooth partition of unity}
\label{sec2}

In this section we adjust the initial sum $W$ to simplify the estimation of type III in \hyperref[sec5]{Section~5}. This technique is also described in~\cite[Section~3]{21}. Suppose that $1 \les a < q \les Q$, $(a,q)=1$. We consider the sum
$$
	W = W(Y) = \sum_{\substack{X \les n < Y \\ n \equiv a\pmod{q}}} \Lambda(n) e(h n^{\alpha}), \qquad X < Y \les 2X.
$$ Let us denote $y = Y/X > 1$. Fix $B_0 > 0$ and choose $\Delta = (\log X)^{-B_0}$.
There exists function $\psi(x) \in C^{\infty}$, such that $\psi(x) = 1$ if $1 \les x \les y$,
$0 \les \psi(x) \les 1$ if $1 - \Delta \les x \les 1$ or $y \les x \les y + \Delta$ and $\psi(x) = 0$ otherwise, and its derivatives satisfy the estimates $\psi^{(j)} (x) \ll_j (\log X)^{j B_0}$. See, for example,~\cite{7} or~\cite{8}. Then $W$ can be rewritten as
\begin{equation} \label{W_expression}
	W = \sum_{\substack{n = 1 \\ n \equiv a\pmod{q}}}^{+\infty} \psi \biggl( \frac{n}{X} \biggr) \Lambda(n) e\bigl(hn^{\alpha}\bigr) + O \biggl( \frac{X (\log X)^{-B_0 + 1}}{q} \biggr). 
\end{equation} By partial summation to prove \hyperref[thm1]{Theorem~1} it is enough to show that the sum in~\eqref{W_expression} is bounded by $X (\log X)^{-A}$. Thus, one can take $B_0 = A + 1$. 

Applying Heath-Brown identity with $k=5$, $V = X^{1/5}$ (\cite[Lemma 1]{14}), we get
$$
	W = \sum_{j=1}^5 (-1)^{j-1} \binom{5}{j} W_j,
$$ where
$$
	W_j = \sum_{\substack{d_1, \ldots, d_{2j} = 1 \\ d_{j+1}, \ldots, d_{2j} \les V \\ d_1 \ldots d_{2j} \equiv a \pmod{q}}} (\log d_1) \mu(d_{j+1}) \ldots \mu(d_{2j}) \psi \biggl( \frac{d_1 \ldots d_{2j}}{X} \biggr) e \bigl( h (d_1 \ldots d_{2j})^{\alpha} \bigr).
$$ The statement of the theorem clearly follows from the estimates $W_j \ll X (\log X)^{-A}$ for each $1 \les j \les 5$.
We only provide the details for $W_5$. The sums $W_1, \ldots, W_4$ can be treated similarly.

We first split the summation over $d_1, \ldots, d_{10}$ to the ``refined'' dyadic intervals following the technique from~\cite{21}. Fix $A_0 > 0$ and $\Theta = 1 + (\log X)^{-A_0}$. Let $\Psi(x)$ be $C^{\infty}$ function supported on $[-\Theta; \Theta]$ such that $\Psi(x) = 1$ on $[-1; 1]$ and $|\Psi^{(j)} (x)| \ll \log^{jA_0} x$ for all $j \ges 0$. For all $x \ges 1$ we have
$$
	1 = \sum_{D \in \mathbf{G}} \Psi_D (x),
$$ where 
$$
	\mathbf{G} = \bigl\{ \Theta^l, \ l \in \mathbb{N} \cup \{0\} \bigr\}, \qquad
	\Psi_D (x) = \Psi \biggl( \frac{x}{D} \biggr) - \Psi \biggl( \frac{\Theta x}{D} \biggr).
$$ Indeed, if $x \ges 1$, then
\begin{multline*}
	\sum_{D \in \mathbf{G}} \Psi_D (x) = \lim_{m \to +\infty} \biggl( \Psi(x) - \Psi(\Theta x) + \Psi\bigl(\frac{x}{\Theta}\bigr) - \Psi(x) + \Psi\bigl(\frac{x}{\Theta^2}\bigr) - \Psi\bigl(\frac{x}{\Theta}\bigr) + \ldots \\ \ldots + \Psi\bigl(\frac{x}{\Theta^m}\bigr) - \Psi\bigl(\frac{x}{\Theta^{m-1}}\bigr) \biggr) = \lim_{m \to +\infty} \biggl( -\Psi(\Theta x) + \Psi\bigl(\frac{x}{\Theta^m}\bigr) \biggr) = -0+1 = 1.
\end{multline*}

The function $\Psi_D$ is supported on $[\Theta^{-1} D; \Theta D]$. Thus,
\begin{multline} \label{W_5}
	W_5 =  \sum_{D_1, \ldots, D_{10} \in \mathbf{G}} \sum_{\substack{d_1, \ldots, d_{10} = 1 \\ d_6, \ldots, d_{10} \les V \\ d_1 \ldots d_{10} \equiv a\pmod{q}}}^{+\infty} \log(d_1) \mu(d_6) \ldots \\ \ldots \mu(d_{10}) \Psi_{D_1} (d_1) \ldots \Psi_{D_{10}} (d_{10}) \psi \biggl( \frac{d_1 \ldots d_{10}}{X} \biggr) e\bigl( h(d_1 \ldots d_{10})^{\alpha} \bigr).
\end{multline} The non-zero contribution to $W_5$ is only coming from the terms satisfying
\begin{equation} \label{d_i_support}
	(1-\Delta) X \les d_1 \ldots d_{10} \les (y+\Delta) X, \qquad \frac{D_i}{\Theta} \les d_i \les D_i \Theta
\end{equation} for $i = 1, \ldots, 10$. From~\eqref{d_i_support} we conclude that the non-zero contribution corresponds
to the tuples $\mathbf{D} = \{D_1, \ldots, D_{10}\}$ satisfying the inequality 
$$
	X_1 \les D_1 \ldots D_{10} \les Y_1, \qquad \text{where} \quad X_1 = (1-\Delta) \Theta^{-10} X, \qquad Y_1 = (y+\Delta) \Theta^{10} X, 
$$ and also satisfying 
$$
	D_i \les V\Theta \qquad \text{for} \quad i = 6, \ldots, 10.
$$ Next, each of the sums $W_1$, $W_2$, $W_3$, $W_4$ and $W_5$ decomposes to the sums of three types depending on the sizes of the variables $d_i$. The decomposition is given by the following lemma:

\begin{lemma}[\textit{\cite[Lemma~3.1]{21}}] \label{lemma1}
	Let $1/10 < \sigma < 1/2$, and let $t_1, \ldots, t_n$ be non-negative real numbers such that $t_1 + \ldots + t_n = 1$. Then at least on of the following three statements holds: \\
	
	{\bfseries Type I:} There is a $t_i$ with $t_i \ges 1/2 + \sigma$;
	
	{\bfseries Type II:} There is a partition $\{1, \ldots, n\} = \mathbf{S} \cap \mathbf{T}$ such that
	$$
		\frac{1}{2} - \sigma < \sum_{i \in \mathbf{S}} t_i \les \sum_{i \in \mathbf{T}} t_i < \frac{1}{2} + \sigma;
	$$
	
	{\bfseries Type III:} There exist distinct $i,j,k$ with $2\sigma \les t_i \les t_j \les t_k \les 1/2 - \sigma$ and
	$$
		t_i + t_j, \ t_i + t_k, \ t_j + t_k \ges \frac{1}{2} + \sigma.
	$$ If $\sigma > 1/6$, then the type III situation is impossible.
\end{lemma} Applying this lemma with $\sigma = 1/10 + \varepsilon_1$, $\varepsilon_1 < 3 \alpha / 5$, to $W_5$ we get
$$
	W_5 \ll |W_I| + |W_{II}| + |W_{III}|,
$$ where the sums correspond to the following cases: \\

{\bfseries Type I sum:} there is one index $1 \les i \les 5$ such that $D_i \ges {X_1}^{3/5 + \varepsilon_1}$.

{\bfseries Type II sum:} there is a partition $\mathbf{S} \cup \mathbf{T} = \{1, \ldots, 10\}$ such that
$$
	{X_1}^{2/5 - \varepsilon_1} < \prod_{i \in \mathbf{S}} D_i < {X_1}^{3/5 + \varepsilon_1}.
$$

{\bfseries Type III sum:} there are three distinct indices $i,j,k \in \{ 1, \ldots, 5 \}$ such that
\begin{gather*}
{X_1}^{1/5 + 2\varepsilon_1} \les D_i \les D_j \les D_k \les {X_1}^{2/5 - \varepsilon_1}, \\
D_i D_j, \ D_i D_k, \ D_j D_k \ges {X_1}^{3/5 + \varepsilon_1}.
\end{gather*}

\begin{remark}
	Note that in the expression analogous to~\eqref{W_5} for $W_1$ and $W_2$ the type III sum is empty.
\end{remark}

\section{The estimation of type I sums}
\label{sec3}

For simplicity we only consider the case $D_1 \ges {X_1}^{3/5 + \varepsilon_1}$. The corresponding sum has the form
\begin{gather*}
	W_I = \sum_{\substack{U \les Y_1 X_1^{-3/5 - \varepsilon_1} \\ U \in \mathbf{G}}} \multsum_{\substack{D_2 \ldots D_{10} = U \\ D_2, \ldots, D_{10} \in \mathbf{G}}} \sum_{\substack{X_1^{3/5 + \varepsilon_1} \les D_1 \les Y_1 U^{-1} \\ D_1 \in \mathbf{G}}} W(\mathbf{D}), \qquad \mathbf{D} = \bigl\{ D_1, \ldots, D_{10} \bigr\}, \\
	W(\mathbf{D}) = \sum_{U \Theta^{-9} \les u \les U \Theta^9} b(u) \sum_{\substack{d_1 = 1 \\ ud_1 \equiv a \pmod{q}}}^{+\infty} f(d_1) e\bigl( h (ud_1)^{\alpha} \bigr),
\end{gather*} where
\begin{gather*}
	b(u) = \sum_{\substack{d_2 \ldots d_{10} = u \\ d_6, \ldots, d_{10} \les V}} \mu(d_6) \ldots \mu(d_{10}) \Psi_{D_2} (d_2) \ldots \Psi_{D_{10}} (d_{10}), \qquad |b(u)| \les \tau_9 (u), \\
	f(d_1) = (\log d_1) \Psi_{D_1} (d_1) \psi \biggl( \frac{ud_1}{X} \biggr).
\end{gather*} Note that the sum over $U \in \mathbf{G}$ contains only $O\bigl((\log X)^{A_0 + 1}\bigr)$ terms. We have
$$
	\bigl|W(\mathbf{D})\bigr| \les || b ||_{\infty} \sum_{U\Theta^{-9} \les u \les U \Theta^9} \biggl| \sum_{\substack{R_1 < d_1 \les R_2 \\ ud_1 \equiv a\pmod{q}}} f(d_1) e \bigl( h (ud_1)^{\alpha} \bigr) \biggr|, 
$$ where 
$$
	|| b ||_{\infty} = \max_{n \les Y_1} |b(n)|, \qquad R_1 = \max \biggl( (1-\Delta)\frac{X}{u}, D_1 \Theta^{-1} \biggr), \qquad R_2 = \min \biggl( (y+\Delta) \frac{X}{u}, D_1 \Theta \biggr). 
$$ By partial summation,
\begin{multline*}
	\bigl|W(\mathbf{D})\bigr| \les || b ||_{\infty} \sum_{U \Theta^{-9} \les u \les U \Theta^{9}} \biggl| f(R_2) \sum_{\substack{R_1 < d_1 \les R_2 \\ ud_1 \equiv a\pmod{q}}} e \bigl( h (ud_1)^{\alpha} \bigr) - \\ \int_{R_1}^{R_2} \biggl( \sum_{\substack{R_1 < d_1 \les v \\ ud_1 \equiv a\pmod{q}}} e \bigl( h(ud_1)^{\alpha} \bigr) \biggr) \frac{df(v)}{dv} dv \biggr|.
\end{multline*} Next,
\begin{multline*}
	\frac{d}{dv} \biggl( \log(v) \Psi_{D_1} (v) \psi \bigl( \frac{uv}{X} \bigr) \biggr) \ll \frac{1}{v} + \frac{\log v}{D_1} (\log X)^{A_0} + \log(v) \frac{u}{X} (\log X)^{B_0} \ll \\ \frac{1}{D_1} (\log X)^{\max(A_0, B_0) + 1},
\end{multline*} and therefore
\begin{multline*}
	\int_{R_1}^{R_2} \biggl( \sum_{\substack{R_1 < d_1 \les v \\ ud_1 \equiv a\pmod{q}}} e \bigl( h(ud_1)^{\alpha} \bigr) \biggr) \frac{df(v)}{dv} dv \ll \\ (\log X)^{\max(A_0, B_0) + 1} \biggl| \sum_{\substack{R_1 < d_1 \les R_3 \\ ud_1 \equiv a\pmod{q}}} e \bigl( h (ud_1)^{\alpha} \bigr) \biggr|,
\end{multline*} where $R_1 < R_3 \les R_2$. Thus, by the triangle inequality,
\begin{equation} \label{Type_I_after_Abel}
	\bigl|W(\mathbf{D})\bigr| \ll || b ||_{\infty} (\log X)^{\max(A_0, B_0) + 1} \sum_{U \Theta^{-9} \les u \les U \Theta^{9}} \biggl| \sum_{\substack{R_1 < d_1 \les R_3 \\ ud_1 \equiv a\pmod{q}}} e \bigl( h (ud_1)^{\alpha} \bigr) \biggr|.
\end{equation} Due to the congruence restriction $ud_1 \equiv a\pmod{q}$ we can assume $(u,q) = 1$ and define $l_1 \equiv au^{\ast} \pmod{q}$, $1 \les l_1 \les q-1$. Setting $d_1 = qr_1 + l_1$, we obtain
$$
	\frac{R_1}{q} \les r_1 + \xi < \frac{R_3}{q}, \qquad \xi = \frac{l_1}{q}.
$$ The inner sum over $d_1$ in~\eqref{Type_I_after_Abel} takes the form
$$
	\sum_{R_1 / q - \xi \les r_1 < R_3 / q - \xi} e \bigl( f_I (r_1) \bigr),
$$ where $f_I (x) = h (uq)^{\alpha} (x + \xi)^{\alpha}$. Then, for $R_1 /q - \xi \les x < R_3 / q - \xi$,
$$
	\bigl|f_I^{''} (x)\bigr| \asymp \frac{h u^{\alpha} q^2}{{R_1}^{2-\alpha}} =: \lambda_2.
$$ By van der Corput second derivative test~\cite[Ch.~1, Theorem 5]{16}, we obtain
$$
	\biggl| \sum_{R_1 / q - \xi \les r_1 < R_3 / q - \xi} e\bigl( f_I (r_1) \bigr) \biggr| \ll \frac{R_3 - R_1}{q} \lambda_2^{1/2} + \lambda_2^{-1/2} = \frac{R_1}{q} \cdot \sqrt{h} \frac{u^{\alpha/2} q}{R_1^{1-\alpha/2}} + \frac{R_1^{1 - \alpha/2}}{\sqrt{h} u^{\alpha/2} q}.
$$ Since $(\log X)^{\max(A_0, B_0) + 1} ||b||_{\infty} \ll_{\delta_1} X^{\delta_1}$ for arbitrary small $\delta_1 > 0$, we get
\begin{multline*}
	\bigl|W(\mathbf{D})\bigr| \ll_{\delta_1} X^{\delta_1} \sum_{U \Theta^{-9} \les u \les U \Theta^{ 9}} \biggl( \frac{R_1}{q} \sqrt{h} \frac{u^{\alpha/2} q}{R_1^{1 - \alpha / 2}} + \frac{R_1^{1 - \alpha/2}}{\sqrt{h} u^{\alpha / 2}q} \biggr) \ll \\ X^{\delta_1} \sqrt{h} R_1^{\alpha / 2} \sum_{U\Theta^9 \les u \les U \Theta^{-9}} u^{\alpha/2} + \frac{X^{\delta_1} R_1^{1 - \alpha/2}}{\sqrt{h} q} \sum_{U\Theta^{-9} \les u \les \Theta U^9} \frac{1}{u^{\alpha/2}} \ll \\ X^{\delta_1} \sqrt{h} D_1^{\alpha / 2} U^{\alpha / 2 + 1} + \frac{X^{\delta_1} D_1^{1 - \alpha / 2} U^{1 - \alpha / 2}}{\sqrt{h}q}.
\end{multline*} Thus,
\begin{multline*}
	W_I \ll \sum_{\substack{U \les Y_1 X_1^{-3/5 - \varepsilon_1} \\ U \in \mathbf{G}}} \multsum_{\substack{D_2 \ldots D_{10} = U \\ D_2, \ldots, D_{10} \in \mathbf{G}}} \sum_{\substack{X_1^{3/5 + \varepsilon_1} \les D_1 \les Y_1 U^{-1} \\ D_1 \in \mathbf{G}}} \biggl( X^{\delta_1} \sqrt{h} D_1^{\alpha / 2} U^{\alpha / 2 + 1} + \\ \frac{X^{\delta_1} D_1^{1 - \alpha / 2} U^{1 - \alpha / 2}}{\sqrt{h}q} \biggr) \ll \sum_{\substack{U \les Y_1 X_1^{-3/5 - \varepsilon_1} \\ U \in \mathbf{G}}} \multsum_{\substack{D_2 \ldots D_{10} = U \\ D_2, \ldots, D_{10} \in \mathbf{G}}} \biggl( \sqrt{h} X^{\delta_1 + \alpha/2} U \log (Y_1 / U)^{A_0 + 1} + \\ \frac{X^{\delta_1 + 1 - \alpha / 2}}{\sqrt{h} q} \log (Y_1 / U)^{A_0 + 1} \biggr).
\end{multline*} Finally, for fixed $U = \Theta^k$, $k \les \log (Y_1 X_1^{-3/5}) / \log \Theta$, using the trivial bound
$$
	\multsum_{\substack{D_2 \ldots D_{10} = U \\ D_2, \ldots, D_{10} \in \mathbf{G}}} 1 = \multsum_{k_2 + \ldots + k_{10} = k} 1 \les k^9 \ll (\log X)^{9(A_0 + 1)},
$$ we get
\begin{multline} \label{contribution_type_I}
W_I \ll \sum_{\substack{U \les Y_1 X_1^{-3/5 - \varepsilon_1} \\ U \in \mathbf{G}}} X^{\delta_1} (\log X)^{10(A_0 + 1)} \biggl( \sqrt{h} U X^{\alpha / 2} + \frac{1}{\sqrt{h}} \frac{X^{1-\alpha / 2}}{q} \biggr) \ll \\ X^{2\delta_1} \biggl( X^{2/5 + \alpha/2 - \varepsilon_1} + \frac{X^{1-\alpha/2}}{q} \biggr).
\end{multline}

\section{The estimation of type II sums}
\label{sec4}

For a fixed partition $\mathbf{S} \cup \mathbf{T} = \{ 1, \ldots, 10 \}$ we use the notation
$$
	m = \prod_{i \in \mathbf{S}} d_i, \qquad n = \prod_{i \in \mathbf{T}} d_i, \qquad M = \prod_{i \in \mathbf{S}} D_i, \qquad N = \prod_{i \in \mathbf{T}} D_i.
$$ Note that $MN \asymp X$. Then type II sum can be written as
$$
	W_{II} = \sum_{\substack{X_1^{2/5 - \varepsilon_1} \les M \les X_1^{3/5 + \varepsilon_1} \\ M \in \mathbf{G}}} \sum_{\substack{X_1 / M \les N \les Y_1 / N \\ N \in \mathbf{G}}} \sum_{\substack{\exists \mathbf{S}, \mathbf{T} : \\ \prod_{i \in \mathbf{S} } D_i = M \\ \prod_{i \in \mathbf{T}} D_i = N}} W(\mathbf{D}),
$$ where
\begin{gather*}
	W(\mathbf{D}) = \sum_{M \Theta^{-|\mathbf{S}|} \les m \les M \Theta^{|\mathbf{S}|}} \gamma(m) \sum_{\substack{n=1 \\ mn \equiv a\pmod{q}}}^{+\infty} \beta(n) \psi \biggl( \frac{mn}{X} \biggr) e\bigl( h(mn)^{\alpha} \bigr), \\
	\gamma(m) = \sum_{\substack{\prod_{i \in \mathbf{S}} d_i = m \\ d_i \les V \ \text{for} \ i \ges 6, i \in \mathbf{S}}} \biggl( \prod_{i \in \mathbf{S}} a_i (d_i) \Psi_{D_i} (d_i) \biggr), \\
	\bigl|\gamma (m)\bigr| \les (\log X) \sum_{\substack{\prod_{i \in \mathbf{S}} d_i = m \\ d_i \les V \ \text{for} \ i \ges 6, i \in \mathbf{S}}} 1 \les (\log X) \tau_{|\mathbf{S}|} (m), \\
	\beta(n) = \sum_{\substack{\prod_{i \in \mathbf{T}} d_i = n \\ d_i \les V \ \text{for} \ i \ges 6, i \in \mathbf{T}}} \biggl( \prod_{i \in \mathbf{T}} a_i (d_i) \Psi_{D_i} (d_i) \biggr), \qquad |\beta(n)| \les (\log X) \tau_{|\mathbf{T}|} (n), \\
	a_1 (d) = \log d, \qquad a_2 (d) = \ldots = a_5 (d) = 1, \qquad a_6 (d) = \ldots = a_{10} (d) = \mu(d).
\end{gather*} By definition of $\beta(n)$ we have
\begin{gather*}
	W(\mathbf{D}) = \sum_{M_1 \les m \les M_2} \gamma(m) \sum_{\substack{N_1 \les n \les N_2 \\ mn \equiv a\pmod{q}}} \beta(n) \psi \biggl( \frac{mn}{X} \biggr) e\bigl( h(mn)^{\alpha} \bigr), \\
	M_1 = M \Theta^{-|\mathbf{S}|}, \qquad M_2 = M \Theta^{|\mathbf{S}|}, \qquad N_1 = N \Theta^{-|\mathbf{T}|}, \qquad N_2 = N \Theta^{|\mathbf{T}|}.
\end{gather*} Cauchy's inequality yields:
$$
	\bigl|W(\mathbf{D})\bigr|^2 \les  \biggl( \sum_{M_1 \les m \les M_2} |\gamma (m)|^2 \biggr) 
	\biggl( \sum_{M_1 \les m \les M_2} \biggl| \sum_{\substack{N_1 \les n \les 
	N_2 \\ mn \equiv a\pmod{q}}} \beta (n) \psi \biggl( \frac{mn}{X} \biggr) e \bigl( h (mn)^{\alpha} \bigr) \biggr|^2 \biggr).		
$$ Next, by Mardzhanishvili's inequality~\cite{20} we get
\begin{equation} \label{Mardz}
	\bigl|W(\mathbf{D})\bigr|^2 \ll M_1 (\log X)^{2 + \kappa}
	\biggl( \sum_{M_1 \les m \les M_2} \biggl| \sum_{\substack{N_1 \les n \les N_2 \\ mn 
	\equiv a\pmod{q}}} \beta (n) \psi \biggl( \frac{mn}{X} \biggr) e \bigl( h (mn)^{\alpha} \bigr) \biggr|^2 \biggr), 
\end{equation} where $\kappa = |\mathbf{S}|^2 - 1$. Rewrite the second factor as follows:
\begin{multline*}
	\sum_{M_1 \les m \les M_2} \sum_{\substack{N_1 \les n_1, n_2 \les 
	N_2 \\ mn_i \equiv a\pmod{q}, i=1,2}} \beta (n_1) \beta (n_2) \psi \biggl( \frac{mn_1}{X} \biggr) \psi \biggl( \frac{mn_2}{X} \biggr) e \bigl( h m^{\alpha} (n_1^{\alpha} 
- 	n_2^{\alpha}) \bigr) = \\ \sum_{M_1 \les m \les M_2} \sum_{\substack{N_1 
	\les n \les N_2 \\ mn \equiv a\pmod{q}}} {\beta^2 (n)} \psi^2 \biggl( \frac{mn}{X} \biggr) + 2 \text{Re} 
	\bigl(S(M, N)\bigr),
\end{multline*} where
$$
	S(M, N) = \sum_{M_1 \les m \les M_2} \sum_{\substack{N_1 \les n_1 < n_2 \les N_2 \\ mn_i 
	\equiv a\pmod{q}, i=1,2}} \beta (n_1) \beta (n_2) \psi \biggl( \frac{mn_1}{X} \biggr) \psi \biggl( \frac{mn_2}{X} \biggr) e \bigl( h m^{\alpha} (n_1^{\alpha} - n_2^{\alpha}) \bigr).
$$ The diagonal term does not exceed
\begin{equation} \label{diagonal}
	\sum_{M_1 \les m \les M_2} \sum_{\substack{N_1 \les n \les N_2 \\ mn \equiv a\pmod{q}}} 
	\beta^2 (n)  \ll (\log X)^2 \frac{MN}{q} (\log X)^{|\mathbf{T}|^2 - 1} \ll \frac{X}{q} (\log X)^{|\mathbf{T}|^2 + 1}.
\end{equation} Setting $m = qr + l$, we get
$$
	R_1 - \eta \les r \les R_2 - \eta, \qquad \eta = \frac{l}{q},
$$ for given $l$, $(l,q) = 1$, $R_1 = M_1 / q$, $R_2 = M_2 / q$. Hence,
\begin{multline*}
	S(M, N) = \sum_{\substack{l=1 \\ (l,q)=1}}^q \sum_{R_1 - \eta \les r \les
	R_2 - \eta} \sum_{\substack{N_1 \les n_1 < n_2 \les 
	N_2 \\ n_1, n_2 \equiv e\pmod{q}}} \beta(n_1) \beta(n_2) \cdot \\ \psi \biggl( \frac{(qr+l)n_1}{X} \biggr) \psi \biggl( \frac{(qr+l)n_2}{X} \biggr) e \bigl( h (n_1^{\alpha} - n_2^{\alpha}) q^{\alpha} (r + \eta)^{\alpha} \bigr),
\end{multline*} where $e = a l^{\ast} \pmod{q}$. Changing the order of summation, we estimate $S(M, N)$ as follows:
\begin{multline*}
	\bigl|S(M, N)\bigr| \les \sum_{\substack{l=1 \\ (l,q)=1}}^q 
	\sum_{\substack{N_1 \les n_1 < n_2 \les N_2 \\ n_1, n_2 \equiv e\pmod{q} }} 
	|\beta (n_1)| |\beta (n_2)| \cdot \\ \biggl| \sum_{R_1 - \eta \les r \les R_2 - \eta} \psi \biggl( \frac{(qr+l)n_1}{X} \biggr) \psi \biggl( \frac{(qr+l)n_2}{X} \biggr) e\bigl( f_{II}(r)\bigr)  \biggr|,
\end{multline*} where $f_{II} (x) = h(n_1^{\alpha} - n_2^{\alpha}) q^{\alpha} (x + \eta)^{\alpha}$. Using the 
conditions $n_1 < n_2$, $n_1 \equiv n_2 \equiv e\pmod{q}$, we set $n_2 = 
n_1 + qs$ with $s \ges 1$. On the other hand, $n_2 \les N_2$ implies $n_1 + 
qs \les N_2$. Hence, $s < (N_2 - N_1) / q = t$, and therefore
\begin{multline*}
	\bigl|S(M, N)\bigr| \ll \sum_{\substack{l=1 \\ (l,q)=1}}^q \sum_{1 \les s < t} 
	\sum_{\substack{N_1 \les n \les N_2 \\ n \equiv e\pmod{q}}} |\beta (n)| |\beta (n+qs)| \cdot \\
	\biggl| \sum_{R_1 - \eta \les r \les R_2 - \eta} \psi \biggl( \frac{(qr+l)n_1}{X} \biggr) \psi \biggl( \frac{(qr+l)n_2}{X} \biggr) e \bigl( f_{II}(r) \bigr) \biggr|.
\end{multline*} By partial summation,
$$
	\bigl|S(M, N)\bigr| \ll (\log X)^{B_0} \sum_{\substack{l=1 \\ (l,q)=1}}^q \sum_{1 \les s < t} 
	\sum_{\substack{N_1 \les n \les N_2 \\ n \equiv e\pmod{q}}} |\beta (n)| |\beta (n+qs)|
	\biggl| \sum_{R_1 - \eta \les r \les R_3 - \eta} e \bigl( f_{II}(r) \bigr) \biggr|,
$$ where $R_1 < R_3 \les R_2$. Next,
$$
	f_{II}^{''} (x) = \frac{\alpha (\alpha - 1) h (n_2^{\alpha} - n_1^{\alpha}) q^{
	\alpha}}{(x + \eta)^{2-\alpha}}, \qquad \text{whence} \quad \bigl|f_{II}^{''} (x)\bigr| \asymp h (n_2^{\alpha} - n_1^{\alpha}) q^{\alpha} \biggl( \frac{q}{M} \biggr)^{2-\alpha}. 
$$ By Lagrange's mean value theorem,
$$
	hq^{\alpha} \bigl( (n + qs)^{\alpha} - n^{\alpha} \bigr) = \alpha hsq^{\alpha+1} 
	(n + qs \theta')^{\alpha - 1} \asymp  hsq^{\alpha+1} N^{\alpha-1} \asymp hsq^{\alpha+1}
	\biggl( \frac{X}{M} \biggr)^{\alpha - 1}, 
$$ where $|\theta'| \les 1$. Hence,
$$
	\bigl|f_{II}^{''}(x)\bigr| \asymp \frac{hsq^2}{X^{1 - \alpha}} \frac{q}{M}.
$$ Applying van der Corput second derivative test~\cite[Ch.~1]{16}, we get
$$
	\sum_{R_1 - \eta \les r \les R_3 - \eta} e \bigl( f_{II} (r) \bigr) \ll \bigl( R_3 - R_1 \bigr) \biggl( \frac{hsq^2}{X^{1-\alpha}} \frac{q}{M} \biggr)^{1/2} + \biggl(\frac{X^{1-\alpha}}{hsq^2} \frac{M}{q} \biggr)^{1/2}.
$$ We have $R_3 - R_1 \ll M/q$. The factor $|\beta (n)| \cdot |\beta (n+qs)|$ is bounded from above by $(X/q)^{\delta_2}$ for arbitrary small $\delta_2 > 0$. The summation over $n \equiv e\pmod{q}$ for $N_1 \les n \les N_2$ contributes the factor of at most $X / (Mq) > 1$ (since $M \ll X^{3/5 + 3\alpha/5}$, $q \les X^{2/5 - 3\alpha / 5}$). Thus,
\begin{multline} \label{non_diagonal_II}
\bigl| S(M, N) \bigr| \ll \\ \biggl(\frac{X}{q} \biggr)^{\delta_2} \sum_{\substack{l=1 \\ (l,q)=1}}^{q-1} \sum_{1 \les s < t} \frac{X}{M q} \biggl( \biggl( \frac{hsqM}{X^{1-\alpha}} \biggr)^{1/2} + \biggl( \frac{M X^{1-\alpha}}{hsq^3} \biggr)^{1/2} \biggr) \ll \\ \biggl( \frac{X}{q} \biggr)^{\delta_2} \biggl( \frac{\sqrt{h} X^{2 + \alpha/2}}{qM^2} + \frac{X^{2 - \alpha/2}}{\sqrt{h} q^2 M} \biggr).
\end{multline} Combining \eqref{Mardz}, \eqref{diagonal} and \eqref{non_diagonal_II}, we get
\begin{multline*}
	\bigl|W(\mathbf{D})\bigr| \ll \\ \sqrt{M} (\log X)^{1 + \kappa/2} \biggl( \frac{X}{q} (\log X)^{|\mathbf{T}|^2 + 1} +  \biggl( \frac{X}{q} \biggr)^{\delta_2} \biggl( \frac{\sqrt{h} X^{2 + \alpha/2}}{qM_1} + \frac{X^{2 - \alpha/2}}{\sqrt{h} q^2 M} \biggr) \biggr)^{1/2} \ll \\ X^{\delta_3} \biggl( \biggl( \frac{XM}{q} \biggr)^{1/2} + \frac{X^{1+\alpha/4}}{(qM)^{1/2}} + \frac{X^{1-\alpha/4}}{q} \biggr),
\end{multline*} where we have used the inequality
$$
	\max \bigl((\log X)^{\kappa / 2 + |\mathbf{T}|^2 + 2}, (X/q)^{\delta_2} \sqrt{h} \bigr) \les X^{\delta_3} 
$$ for some $\delta_3 \ges \delta_2$. For fixed $M = \Theta^k$ and $N = \Theta^l$ with $k+l = 10$ the number of corresponding tuples $\mathbf{S}$ and $\mathbf{T}$ does not exceed
$$
	\sum_{\substack{i+j=10 \\ i+j \ges 1}} k^i l^j \ll \biggl( \frac{\log X}{\log \Theta} \biggr)^{10} \ll (\log X)^{10(A_0 + 1)}.
$$ Thus,
\begin{multline} \label{contribution_type_II}
	W_{II} \ll \sum_{\substack{X_1^{2/5 - \varepsilon_1} \les M \les X_1^{3/5 + \varepsilon_1} \\ M \in \mathbf{G}}} \sum_{\substack{X_1 / M \les N \les Y_1 / N \\ N \in \mathbf{G}}} \sum_{\substack{\exists \mathbf{S}, \mathbf{T} : \\ \prod_{i \in \mathbf{S} } D_i = M \\ \prod_{i \in \mathbf{T}} D_i = N}} \bigl|W(\mathbf{D})\bigr| \ll \\
	(\log X)^{10(A_0 + 1)} \sum_{\substack{X_1^{2/5 - \varepsilon_1} \les M \les X_1^{3/5 + \varepsilon_1} \\ M \in \mathbf{G}}} \sum_{\substack{X_1 / M \les N \les Y_1 / N \\ N \in \mathbf{G}}} X^{\delta_3} \biggl( \biggl( \frac{XM}{q} \biggr)^{1/2} + \frac{X^{1+\alpha/4}}{(qM)^{1/2}} + \frac{X^{1-\alpha/4}}{q} \biggr) \ll \\
	X^{2 \delta_3} \biggl( \frac{X^{4/5+\varepsilon_1 / 2}}{\sqrt{q}} + \frac{X^{4/5 + \alpha / 4 + \varepsilon_1 / 2}}{\sqrt{q}} + \frac{X^{1 - \alpha / 4}}{q} \biggr) \ll \frac{X^{4/5 + \alpha / 4 + \varepsilon_1 / 2 + 2\delta_3}}{\sqrt{q}} + \frac{X^{1 - \alpha / 4 + 2 \delta_3}}{q}.
\end{multline}

\section{The estimation of type III sums}
\label{sec5}

We apply the method of stationery phase to treat the type III sum. To deal with the oscillatory integrals arising after the Poisson summation we use two auxiliary lemmas given below:

\begin{lemma}[\textit{Lemma 8.1, \cite{3}}] \label{lemma2}
	Let $Y_I \ges 1$, $X_I, Q_I, V_I, R_I > 0$, $w(t)$ is a smooth function supported on some finite interval $\mathbb{J} \subset \mathbb{R}$ such that 
	$$
		w^{(j)}(t) \ll_j X_I V_I^{-j}
	$$ for all $j \ges 0$. Suppose that $g(t)$ is a smooth function such that $|g'(t)| \ges R_I$, $g^{(j)} (t) \ll_j Y_I Q_I^{-j}$ for $j \ges 2$, $t \in \mathbb{J}$. Then the integral $I$ defined by
	$$
		I = \int_{-\infty}^{+\infty} w(t) e\bigl( g (t)\bigr) dt
	$$ satisfies
	$$
		I \ll_{A_I} |\mathbb{J}| X_I \bigl( (Q_I R_I / \sqrt{Y_I})^{-A_I} + (R_I V_I)^{-A_I} \bigr)
	$$ with any fixed real $A_I > 0$.
\end{lemma} This result gives a non-trivial upper bound for the integral $I$ in the case if $R_I V_I$ and $Q_I R_I Y_I^{-1/2}$ are much bigger than 1. \\

\begin{lemma}[\textit{Proposition 8.2, \cite{3}}] \label{lemma3}
	Let $0 < \delta_I < 1/10$, $X_I, Y_I, V_I, {\tilde V_I}$, $Q_I > 0$, $Z_I = Q_I + X_I + Y_I + {\tilde V_I} + 1$, and assume that $Y_I \ges Z_I^{3\delta_I}$,
	$$
		{\tilde V_I} \ges V_I \ges \frac{Q_I Z_I^{\delta_I / 2}}{Y_I^{1/2}}.
	$$ Suppose that $w(t)$ is a smooth function supported on an interval $\mathbb{J}$ of length ${\tilde V_I}$ satisfying
	$$
		w^{(j)}(t) \ll_j X_I V_I^{-j}
	$$ for all $j \ges 0$. Suppose that $g(t)$ is a smooth function such that there is unique point $t_0 \in \mathbb{J}$ such that $g'(t_0) = 0$. Further, $g(t)$ satisfies the estimates $g''(t) < 0$, $g''(t) \gg Y_I Q_I^{-2}$, $g^{(j)}(t) \ll_j Y_I Q_I^{-j}$, for all $j \ges 1$, $t \in \mathbb{J}$. Then the integral
	$$
		I = \int_{-\infty}^{+\infty} w(t) e\bigl(g(t)\bigr) dt
	$$ has an asymptotic expansion of the form
	\begin{gather*}
		I = \frac{e\bigl(g(t_0)\bigr)}{|g''(t_0)|^{1/2}} \sum_{0 \les n \les 3\delta_I^{-1} A_I} p_n (t_0) + O_{A_I, \delta_I} (Z_I^{-A_I}), \\
		p_n (t_0) = \frac{\sqrt{2\pi} e^{-\pi i /4}}{n!} \frac{(2i)^{-n}}{|g''(t_0)|^n} G^{(2n)} (t_0),
	\end{gather*} where $A_I > 0$ is arbitrary, and
	$$
		G(t) = w(t) e\bigl(H(t)\bigr), \qquad H(t) = g(t) - g(t_0) - \frac{1}{2} g''(t_0) (t-t_0)^2.
	$$ 
\end{lemma} 

Later in this section we will use the following notation:
\begin{gather*}
	\beta = \frac{2-\alpha}{1-\alpha}, \qquad \gamma = \frac{\alpha}{1-\alpha}, \qquad \delta = \frac{1}{1-\alpha}, \\
	\xi = \frac{1}{1-\gamma} = \frac{1-\alpha}{1-2\alpha}, \qquad \eta = \frac{\alpha}{1-2\alpha}, \qquad \omega = \xi (2-\gamma) = \frac{2-3\alpha}{1-2\alpha}.
\end{gather*} Let us denote as $M,N,K$ the three indices from $\{ D_1, \ldots, D_5 \}$ satisfying type III conditions, as $m,n,k$ the corresponding indices from $\{ d_1, \ldots, d_5 \}$, as $i_1, i_2, i_3$ the corresponding indices from $\{ 1, \ldots, 5\}$, and let $\mathbf{I}$ be the set of all remaining indices $\{ 1, \ldots, 10\}$\textbackslash$\{ i_1,i_2,i_3 \}$. Also let
$$
	U = \prod_{i \in \mathbf{I}} D_i, \qquad u = \prod_{i \in \mathbf{I}} d_i.
$$ We get the sum of the form
\begin{multline} \label{type_III}
	W_{III} = \mathop{{\sum}^{'}}_{M,N,K \in \mathbf{G}} \mathop{{\sum}^{'}}_{U \in \mathbf{G}} \sum_{U\Theta^{-7} \les u \les U\Theta^7} F(U, u) \sum_{\substack{m,n,k = 1 \\ umnk \equiv a\pmod{q}}}^{+\infty} f_1 (m) f_2 (n) f_3 (k)  \cdot \\ \Psi_M (m) \Psi_N (n) \Psi_K (k) \psi \biggl( \frac{umnk}{X} \biggr) e\bigl( h (umnk)^{\alpha} \bigr),
\end{multline} where
\begin{gather*}
	F(U, u) = \biggl( \multsum_{\prod_{i \in \mathbf{I}} D_i = U} \biggr) \biggl( \multsum_{\prod_{i \in \mathbf{I}} d_i = u} \biggr) \biggl( \prod_{i \in \mathbf{I}} a_i (d_i) \Psi_{D_i} (d_i) \biggr), \\
a_1 (d) = \log d, \qquad a_2 (d) = \ldots = a_5 (d) = 1, \qquad a_6 (d) = \ldots = a_{10} (d) = \mu(d),
\end{gather*} $f_i (x)$ are smooth functions such that $f_i (x) \equiv 0$ if $x \les 0$ and $f_i (x) = 1$ or $f_i (x) = \log x$ for $x \ges 1$, $\sum'$ denotes the summation over $M,N,K,U \in \mathbf{G}$ satisfying the type III conditions. Without loss of generality we can assume $M \les N \les K$. Then rewrite~\eqref{type_III} in the following way:
\begin{multline} \label{W_III_new_form}
	W_{III} = \\ \sum_{\substack{U \les X^{1/10 - 3\varepsilon_1 / 2} \\ U \in \mathbf{G}}} \sum_{U \Theta^{-7} \les u \les U \Theta^7} F(U,u) \sum_{\substack{M_1 \les M \les M_2 \\ M \in \mathbf{G}}} \sum_{\substack{N_1 \les N \les N_2 \\ N \in \mathbf{G}}} \sum_{\substack{K_1 \les K \les K_2 \\ K \in \mathbf{G}}} W(M,N,K),
\end{multline} where
\begin{gather*}
	M_1 = X_1^{1/5 + 2\varepsilon_1}, \qquad M_2 = (Y_1 U^{-1} )^{1/3}, \\
	N_1 = \max \bigl( M, X^{3/5 + \varepsilon_1} M^{-1} \bigr), \qquad N_2 = \min \biggl(X^{2/5 - \varepsilon_1}, \biggl( \frac{Y}{MU} \biggr)^{1/2} \biggr), \\
	K_1 = N, \qquad K_2 = \min \biggl( X^{2/5 - \varepsilon_1}, \frac{Y}{UMN} \biggr),
\end{gather*} and
\begin{multline*}
W(M,N,K) = \\ \sum_{\substack{m,n,k = 1 \\ mnku \equiv a\pmod{q}}}^{+\infty} f_1 (m) f_2 (n) f_3 (k) \Psi_M (m) \Psi_{N} (n) \Psi_K(k) \psi \biggl(\frac{umnk}{X} \biggr) e \bigl( h (umnk)^{\alpha} \bigr).
\end{multline*} Note that $f_i (.) \Psi_D (.) \psi(.)$ is smooth on $(0; +\infty)$, so one can apply Poisson summation to any of the sums over $n,m,k$. We also note that the number of terms in each sum over $U, M, N, K \in \mathbf{G}$ in~\eqref{W_III_new_form} is $O\bigl((\log X)^{A_0 + 1} \bigr)$ and $|F(U, u)|$ can be bounded as follows:
\begin{multline*}
	\bigl|F(U,u)\bigr| \ll (\log X) \tau_7 (u) \cdot \# \biggl\{ (e_1, \ldots, e_7) \in \mathbb{Z}_{\ges 0}^7 : e_1 + \ldots + e_7 = \frac{\log U}{\log \Theta} \biggr\} \ll \\  (\log X) \tau_7 (u) \biggl( \frac{\log U}{\log \Theta} \biggr)^6 \ll \tau_7 (u) (\log X)^{6(A_0 + 1) + 1}.
\end{multline*}

\subsection*{First iteration of Poisson summation}

We first apply Poisson summation to the longest sum over $k$. By orthogonality of characters,
$$
	W(M,N,K) = \frac{1}{\varphi(q)} \sum_{\chi \ \text{mod} \ q} \chi(ua^{\ast}) \sum_{m=1}^{+\infty} \chi(m) f_1 (m) \Psi_M (m) \sum_{n=1}^{+\infty} \chi(n) f_2 (n) \Psi_N (n) W_{m,n,\chi},
$$ where
$$
	W_{m,n,\chi} = \sum_{k=1}^{+\infty} \chi(k) f_3 (k) \Psi_K (k) \psi \biggl( \frac{umnk}{X} \biggr) e \bigl( h(umnk)^{\alpha} \bigr).
$$ To remove the factor $\chi(k)$ in the last sum we substitute $k=qr+l$:
$$
	W_{m,n,\chi} = \sum_{l=1}^{q-1} \chi(l) \sum_{r=-\infty}^{+\infty} f_3 (qr+l) \Psi_K (qr+l) \psi \biggl( \frac{umn (qr+l)}{X} \biggr) e \bigl( h (umn(qr+l))^{\alpha}\bigr).
$$ The function $(qr+l)^{\alpha}$ is extended by zero for $r < -l/q$. By Poisson summation,
\begin{multline*}
	W_{m,n,\chi} = \\ \sum_{l=1}^{q-1} \chi(l) \sum_{s=-\infty}^{+\infty} \int_{-\infty}^{+\infty} f_3 (qv+l) \Psi_K (qv+l) \psi \biggl( \frac{umn(qv+l)}{X} \biggr) e \bigl( h(umn(qv+l))^{\alpha} \bigr) e\bigl( -vs\bigr) dv.
\end{multline*} We can reduce the range of integration to $(-l/q; +\infty)$ due to the fact that $f_3 (x) = 0$ for $x \les 0$. Then substituting
$$
	t = \frac{umn(qv+l)}{X}
$$ we get
$$
	W_{n,m,\chi} = \frac{X}{qumn} \sum_{s=-\infty}^{+\infty} \tau(\chi; s) I_{m,n} (s),
$$ where
\begin{gather*}
	I_{m,n} (s) = \int_0^{+\infty} f_3 \biggl( \frac{Xt}{umn} \biggr) \Psi_K \biggl( \frac{Xt}{umn} \biggr) \psi(t) e \biggl( h(Xt)^{\alpha} - \frac{Xst}{qumn} \biggr) dt, \\
	\tau(\chi; s) = \sum_{l=1}^{q-1} \chi(l) e\biggl( \frac{sl}{q} \biggr) \qquad \text{is a Gauss sum.}
\end{gather*} Next, we verify the conditions of \hyperref[lemma2]{Lemma~2} and \hyperref[lemma3]{Lemma~3}. Let
\begin{gather*}
	w(t) = f_3 \biggl( \frac{Xt}{umn} \biggr) \Psi_K\biggl( \frac{Xt}{umn} \biggr) \psi(t), \\
	g_n (t) = 
	\begin{cases}
	\displaystyle h(Xt)^{\alpha} - \frac{Xst}{qumn}, 
	&\text{if }
	\displaystyle 1 - \Delta \les t \les y + \Delta,
	\\[2mm]
	\displaystyle 0
	&\text{if }
	\displaystyle t \les 1-2\Delta \quad \text{or }
	\displaystyle t \ges y + 2\Delta,
	\end{cases}
\end{gather*} and extend $g_n (t)$ to a smooth function on $[1-2\Delta, 1-\Delta]$ and $[y+\Delta, y+2\Delta]$. We now evaluate the derivatives. First, if $1 - \Delta \les t \les y + \Delta$ and $j \ges 2$, then we have
$$
	g_n^{(j)} (t) = \frac{(\alpha)_j hX^{\alpha}}{t^{j-\alpha}}, \qquad (\alpha)_j = \prod_{i=1}^j (\alpha - i + 1), \qquad \bigl|g_n^{(j)}(t)\bigr| \asymp_{\alpha, j} hX^{\alpha}.
$$ Thus, one can take $Y_I = hX^{\alpha}$, $Q_I = 1$. Now let us estimate $w^{(j)}(t)$ on $\mathbb{J}$. We have
$$
	\frac{d^j w(t)}{dt^j} = \sum_{j_1 + j_2 + j_3 = j} \binom{j}{j_1, j_2, j_3} \frac{d^{j_1} f_3}{dt^{j_1}} \biggl( \frac{Xt}{umn} \biggr) \frac{d^{j_2} \Psi_K}{dt^{j_2}} \biggl( \frac{Xt}{umn} \biggr) \frac{d^{j_3} \psi(t)}{dt^{j_3}}.
$$ Next,
\begin{gather*}
	\frac{d^{j_1}f_3}{dt^{j_1}} \biggl( \frac{Xt}{umn} \biggr) \ll \log X, \\
	\frac{d^{j_2}}{dt^{j_2}} \Psi_K \biggl( \frac{Xt}{umn} \biggr) \ll \biggl( \frac{X}{Kumn} \biggr)^{j_2} (\log X)^{j_2 A_0} \ll (\log X)^{j_2 A_0}, \\
	\frac{d^{j_3} \psi(t)}{dt^{j_3}} \ll (\log X)^{j_3 B_0}.
\end{gather*} Thus, we find
\begin{multline} \label{bound_derivative_w}
w^{(j)} (t) \ll (\log X) \sum_{j_1 + j_2 + j_3 = j} \binom{j}{j_1, j_2, j_3} (\log X)^{j_2 A_0} (\log X)^{j_3 B_0} \ll \\ (\log X) \bigl( 1 + (\log X)^{A_0} + (\log X)^{B_0} \bigr)^j \ll (\log X)^{C_0 j + 1}, 
\end{multline} where $C_0 = \max(A_0, B_0)$. So one can take $X_I = \log X$, $V_I = (\log X)^{-C_0}$. From $Z_I = Q_I + X_I + Y_I + \tilde V_I + 1$ we get $Z_I \asymp Y_I \asymp hX^{\alpha}$, whence for any fixed $\delta_I$, $0 < \delta_I < 1/10$,
we have
$$
	V_I = (\log X)^{-C_0} \ges \frac{Q_I Z_I^{\delta_I / 2}}{\sqrt{Y_I}} \asymp (h X^{\alpha})^{-1/2 + \delta_I / 2}.
$$ Now set
$$
	T_1 = \frac{1}{4} \frac{\alpha h u mNq}{X^{1-\alpha}}, \qquad T_2 = 4 \frac{\alpha h u mNq}{X^{1-\alpha}}
$$ and split the sum $W_{m,n,\chi}$ in the following way:
\begin{multline*}
W_{m,n,\chi} = \frac{X}{qumn} \biggl\{ \sum_{T_1 \les s \les T_2} + \sum_{|s| > T_2} + \sum_{-T_2 \les s < T_1} \biggr\} \tau(\chi; s) I_{m,n} (s) =: \\ \frac{X}{qumn} (S_1 + S_2 + S_3).
\end{multline*} For $S_2$ and $S_3$ we apply \hyperref[lemma2]{Lemma~2} to estimate  $I_{m,n} (s)$; for $S_1$ we compute $I_{m,n} (s)$ asymptotically using \hyperref[lemma3]{Lemma~3}. We have
$$
	g_n'(t) = \alpha h X^{\alpha} t^{\alpha - 1} - \frac{Xs}{qumn}.
$$ If $|s| > T_2$, then
$$
	|g_n'(t)| \ges \frac{X|s|}{qumn} \biggl( 1 - \frac{\alpha h t^{\alpha-1} qumn}{X^{1-\alpha} T_2} \biggr) \ges \frac{X|s|}{2qumn}.
$$ If $-T_2 \les s \les 0$, then
$$
	g_n'(t) = \alpha h X^{\alpha} t^{\alpha-1} + \frac{X|s|}{qumn} \ges \alpha h X^{\alpha} t^{\alpha-1} \ges \frac{\alpha}{3} h X^{\alpha}.
$$ Finally, if $1 \les s < T_1$, then
$$
	g_n'(t) \ges \alpha h X^{\alpha} t^{\alpha-1} \biggl( 1 - \frac{XT_1}{qumn} \frac{t^{1-\alpha}}{\alpha h X^{\alpha}} \biggr) \alpha h X^{\alpha} t^{\alpha-1} \biggl( 1- \frac{5}{8} \biggr) \ges \frac{\alpha}{6} h X^{\alpha}.
$$ Thus, one can choose
$$
	R_I = 
	\begin{cases}
	\displaystyle \frac{X|s|}{2qumn}
	&\text{if }\displaystyle |s| > T_2,
	\\[3mm]
	\displaystyle \frac{\alpha}{6} h X^{\alpha} 
	&\text{if } \displaystyle -T_2 \les s < T_1.
	\end{cases}
$$ In the case $|s| > T_2$ we set
$$
	\Delta_1 = \frac{Q_I R_I}{\sqrt{Y_I}}, \qquad \Delta_2 = R_I V_I,
$$ and get
\begin{gather*}
	\Delta_1 = \frac{X|s|}{2qumn} \frac{1}{\sqrt{hX^{\alpha}}} \ges \frac{X^{1-\alpha/2} T_2}{2\sqrt{h} qumn} \ges X^{\alpha/2} \cdot 2\alpha \sqrt{h} \frac{N}{n} \ges \alpha X^{\alpha/2}, \\
	\Delta_2 = \frac{X|s|}{2qumn} (\log X)^{-C_0} \ges \frac{XT_2 (\log X)^{-C_0}}{2qumn} \ges X^{\alpha} \cdot \frac{2\alpha N}{n} (\log X)^{-C_0} \ges X^{\alpha/2}.
\end{gather*} If $-T_2 \les s < T_1$, then
$$ 
	\Delta_1 = \frac{\alpha h}{6} \frac{X^{\alpha}}{\sqrt{hX^{\alpha}}} \ges \frac{\alpha}{6} X^{\alpha/2}, \qquad
	\Delta_2 = \frac{\alpha h}{6} X^{\alpha} (\log X)^{-C_0} \ges X^{\alpha/2}.
$$ Thus, by \hyperref[lemma2]{Lemma~2},
\begin{multline*}
	I_{m,n} (s) \ll_{\alpha} (\log X) \biggl\{ \biggl( \frac{X|s|}{2qumn} \frac{1}{\sqrt{hX^{\alpha}}} \biggr)^{-A_I} + \biggl( \frac{X|s|}{2qumn} \frac{1}{(\log X)^{C_0}} \biggr)^{-A_I} \biggr\} \ll_{\alpha} \\ 
	(\log X) \biggl( \frac{X|s|}{2qumn} \frac{1}{\sqrt{h X^{\alpha}}} \biggr)^{-A_I} \ll_{\alpha} (\log X) \biggl( \frac{2qumn \sqrt{h}}{X^{1-\alpha/2} |s|} \biggr)^{A_I}
\end{multline*} for $|s| \ges T_2$, and
\begin{multline*}
	I_{m,n} (s) \ll_{\alpha} (\log X) \biggl\{ \biggl( \frac{\alpha}{6} \sqrt{h} X^{\alpha/2} \biggr)^{-A_I} + \biggl( \frac{\alpha h}{6} X^{\alpha} (\log X)^{-C_0} \biggr)^{-A_I} \biggr\} \ll_{\alpha} \\
	(\log X) X^{-\alpha A_I / 2}
\end{multline*} if $-T_2 \les s \les T_1$. Choose $A_I = 2D_0 + 1$, where $D_0 = D_0 (\alpha) > 1$ is large enough. Going back to 	$S_2$ and $S_3$, we get
\begin{gather*}
	S_2 \ll \sum_{|s| > T_2} (\log X) \biggl( \frac{2qumn \sqrt{h}}{X^{1-\alpha/2}} \biggr)^{A_I} \frac{1}{|s|^{A_I}} \ll \frac{umNq}{X^{1-\alpha}} X^{-\alpha D_0}, \\
	S_3 \ll (T_1 + T_2 + 1) X^{-\alpha D_0 - \alpha/2} (\log X) \ll \biggl( \frac{umNq}{X^{1-\alpha}} + 1 \biggr) X^{-\alpha D_0}.
\end{gather*} From $qumn \gg XK^{-1}$ we find
\begin{multline*}
	W_{m,n,\chi} = \frac{X}{qumn} S_1 + O \biggl( \frac{X}{qumn} \biggl( \frac{umNq}{X^{1-\alpha}} + 1\biggr) X^{-\alpha D_0} \biggr) = \\
	\frac{X}{qumn} S_1 + O \bigl( X^{-\alpha (D_0 + 1)} + K X^{-\alpha D_0} \bigr) = \frac{X}{qumn} S_1 + O\bigl( KX^{-\alpha D_0}  \bigr).
\end{multline*}  

Now we compute $S_1$. Choosing $\delta_I = 1/20$, $A_I = D_0$, we apply \hyperref[lemma3]{Lemma~3} to $I(s)$ when $T_1 < s \les T_2$. Let $g_n'(t_0) = 0$. Then
$$
	t_0 = \frac{1}{X} \biggl( \frac{\alpha hqumn}{s} \biggr)^{1 / (1-\alpha)}.
$$ Notice that for any $T_1 \les s \les T_2$ the point $t_0$ lies in $\mathbb{J} = \bigl[10^{-1}; 10\bigr]$. Thus,
$$
	I_{m,n} (s) = e \biggl( g_n(t_0) - \frac{1}{8} \biggr) \sum_{0 \les \nu \les \nu_1} \frac{\sqrt{2\pi}}{\nu!} \frac{(2i)^{-\nu}}{|g_n''(t_0)|^{\nu+1/2}} \frac{d^{2\nu} G_n(t)}{dt^{2\nu}} \biggl|_{t=t_0} + O\bigl(X^{-\alpha D_0}\bigr),
$$ where
$$
	G_n(t) = w(t) e\bigl(H_n(t)\bigr), \qquad \nu_1 = 60D_0, \qquad H_n(t) = g_n(t) - g_n(t_0) - \frac{1}{2} g_n''(t_0) (t-t_0)^2.
$$ One can easily verify the identities
$$
	g_n(t_0) = (1-\alpha) (\alpha^{\alpha} h)^{\delta} \biggl( \frac{qumn}{s} \biggr)^{\gamma}, \qquad \bigl|g_n''(t_0)\bigr| = \alpha (1-\alpha) hX^2 \biggl( \frac{s}{\alpha h qumn} \biggr)^{\beta},
$$ where $\gamma = \alpha / (1-\alpha), \delta=1 / (1-\alpha), \beta = (2-\alpha)/ (1-\alpha)$. Then, if $1-\Delta \les t_0 \les y+\Delta$, we get
\begin{multline} \label{expression}
	I_{m,n} (s) = e \biggl( (1-\alpha) (\alpha^{\alpha} h)^{\delta} \biggl( \frac{qumn}{s} \biggr)^{\delta} \biggr) \sum_{0 \les \nu \les \nu_1} \frac{c_{\nu} (\alpha)}{(hX^2)^{\nu+1/2}} \biggl( \frac{hqumn}{s} \biggr)^{\beta (\nu + 1/2)} \cdot \\ \frac{d^{2\nu} G_n(t)}{dt^{2\nu}} \biggl|_{t=t_0} + O\bigl(X^{-\alpha D_0}\bigr),
\end{multline} with
$$
	c_{\nu} (\alpha) = \frac{\sqrt{2\pi}}{\nu!} \frac{(2i)^{-\nu} e^{-\pi i /4}}{\bigl(\alpha (1-\alpha)^{\nu+1/2}\bigr)} \alpha^{\beta(\nu+1/2)}.
$$ Notice that~\eqref{expression} remains valid if $t_0 \notin \bigl[1-\Delta; y+\Delta\bigr]$ since $w(t) \equiv 0$, $G_n(t) \equiv 0$ for $t$ close to $t_0$. 

Going back to the sum $W_{m,n,\chi}$, we have
\begin{multline*}
	W_{m,n,\chi} = \frac{X}{qumn} S_1 + O\bigl(KX^{- \alpha D_0}\bigr) = \frac{X}{qumn} \sum_{T_1 < s \les T_2} \tau(\chi; s) \biggl\{ e \biggl( (1-\alpha) (\alpha^{\alpha} h)^{\delta} \biggl( \frac{qumn}{s} \biggr)^{\gamma} \biggr) \cdot \\ \sum_{0 \les \nu \les \nu_1} \frac{c_{\nu} (\alpha)}{(hX^2)^{\nu + 1/2}} \biggl( \frac{hqumn}{s} \biggr)^{\beta(\nu+1/2)} \frac{d^{2\nu} G_n(t)}{dt^{2\nu}} \biggl|_{t=t_0} + O\bigl(X^{-\alpha D_0}\bigr) \biggr\} + O\bigl(KX^{-\alpha D_0}\bigr). 
\end{multline*} The contribution from the error terms can be made arbitrarily small with the appropriate choice of $D_0$. The main term takes the form
\begin{multline*}
	\frac{1}{\varphi(q)} \sum_{\chi \ \text{mod} \ q} \chi(u a^{\ast}) \sum_{m=1}^{+\infty} \chi(m) f_1 (m) \Psi_M(m) \sum_{n=1}^{+\infty} \chi(n) f_2 (n) \Psi_N(n) \cdot \\ \frac{X}{qumn} \sum_{0 \les \nu \les \nu_1} \frac{c_{\nu} (\alpha)}{(hX^2)^{\nu + 1/2}} \sum_{T_1 \les s \les T_2} \tau(\chi; s) \biggl( \frac{hqumn}{s} \biggr)^{\beta(\nu+1/2)} \frac{d^{2\nu} G_n(t)}{dt^{2\nu}} \biggl|_{t=t_0} \cdot \\  e \biggl\{ (1-\alpha) (\alpha^{\alpha} h)^{\delta} \biggl( \frac{qumn}{s} \biggr)^{\gamma} \biggr\}.
\end{multline*} We also note that for the small values of $q$ it is possible to get $T_2 < 1$. This case is not a problem since the sum $S_1$ is empty and the only contribution to the upper bound is coming from~\hyperref[lemma2]{Lemma~2}.

\subsection*{Second iteration of Poisson summation}

We have:
\begin{multline} \label{V_chi_intro}
	W(M,N,K) = \\ \frac{(qu)^{-1}}{\varphi(q)} \sum_{\chi \ \text{mod} \ q} \chi(u a^{\ast}) \sum_{m=1}^{+\infty} \chi(m) \frac{f_1 (m)}{m} \Psi_M (m) \sum_{T_1 < s < T_2} \tau(\chi; s) V_{\chi, m, s} + O\bigl(X^{-\alpha D_0 / 2}\bigr),
\end{multline} where
\begin{multline*}
	V_{\chi, m, s} = \sum_{n=1}^{+\infty} \chi(n) \frac{f_2 (n)}{n} \Psi_N (n) \sum_{0 \les \nu \les \nu_1} \frac{c_{\nu} (\alpha)}{X^{2\nu}} h^{\delta (\nu + 1/2)} \biggl( \frac{qumn}{s} \biggr)^{\beta (\nu + 1/2)} \cdot \\ \frac{d^{2\nu} G_n(t)}{dt^{2\nu}} \biggl|_{t=t_0} e \biggl\{ (1-\alpha) (\alpha^{\alpha} h)^{\delta} \biggl( \frac{qumn}{s} \biggr)^{\gamma} \biggr\}.
\end{multline*} Setting $n = q\rho + \lambda$, we get
\begin{multline*}
	V_{\chi, m, s}  = \sum_{0 \les \nu \les \nu_1} \frac{c_{\nu} (\alpha)}{X^{2\nu}} h^{\delta (\nu + 1/2)} \biggl( \frac{qum}{s} \biggr)^{\beta (\nu + 1/2)} \sum_{\lambda=1}^{q-1} \chi(\lambda) \sum_{\rho = -\infty}^{+\infty} f_2 (q\rho + \lambda) \Psi_N(q\rho + \lambda) \cdot \\  (q\rho + \lambda)^{\beta(\nu + 1/2) - 1} \frac{d^{2\nu} G_{q\rho + \lambda}(t)}{dt^{2\nu}} \biggl|_{t=t_0} e \biggl\{ (1-\alpha) (\alpha^{\alpha} h)^{\delta} \biggl( \frac{qum (q\rho + \lambda)}{s} \biggr)^{\gamma} \biggr\}.
\end{multline*} Applying Poisson summation again, we obtain
\begin{multline} \label{for_integral}
	V_{\chi, m, s} = \sum_{0 \les \nu \les \nu_1} \frac{c_{\nu} (\alpha)}{X^{2\nu}} h^{\delta (\nu + 1/2)} \biggl( \frac{qum}{s} \biggr)^{\beta (\nu + 1/2)} \cdot \\ \sum_{\lambda=1}^{q-1} \chi(\lambda) \sum_{\sigma = -\infty}^{+\infty} \int_{-\infty}^{+\infty} f_2 (qv+\lambda) \Psi_N (qv+\lambda) (qv+\lambda)^{\beta(\nu + 1/2) - 1} \cdot \\ \frac{d^{2\nu} G_{q\rho + \lambda}(t)}{dt^{2\nu}} \biggl|_{t=t_0} e \biggl\{ (1-\alpha) (\alpha^{\alpha} h)^{\delta} \biggl( \frac{qum (qv+\lambda)}{s} \biggr)^{\gamma} - \sigma v \biggr\} dv.
\end{multline} Next, we substitute
$$
	\tau = \frac{\alpha h qum (qv + \lambda)}{sX^{1-\alpha}}.
$$ This implies
\begin{gather*}
	t_0 = \frac{1}{X} \biggl( \frac{\alpha h q um (qv+\lambda)}{s} \biggr)^{\delta} = \tau^{\delta}, \qquad
	\biggl( \frac{qum (qv+\lambda)}{s} \biggr)^{\gamma} = \biggl( \frac{X^{1-\alpha} \tau}{\alpha h} \biggr)^{\gamma}, \\
	(\alpha^{\alpha} h)^{\delta} \biggl( \frac{qum (qv+\lambda)}{s} \biggr)^{\gamma} = h X^{\alpha} \tau^{\gamma}.  
\end{gather*} For convenience we will further use a slightly different notation for functions $G_n, H_n$ and $g_n$: $G_n (t) = G(t, \tau), H_n (t) = H(t, \tau), g_n(t) = g(t, \tau)$. The integral in~\eqref{for_integral} takes the form
\begin{multline*}
	\frac{1}{q} \biggl( \frac{X^{1-\alpha} s}{\alpha h qum} \biggr)^{\beta(\nu+1/2)} e\biggl( \frac{\lambda \sigma}{q} \biggr) \int_0^{+\infty} f_2 \biggl( \frac{X^{1-\alpha} s\tau}{\alpha h qum} \biggr) \Psi_N \biggl( \frac{X^{1-\alpha} s\tau}{\alpha h qum} \biggr) \tau^{\beta (\nu + 1/2) - 1} \cdot \\ \frac{d^{2\nu} G(t, \tau)}{dt^{2\nu}} \biggl|_{t=\tau^{\delta}} e \biggl\{ (1-\alpha) hX^{\alpha} \tau^{\gamma} - \frac{X^{1-\alpha} s\sigma \tau}{\alpha h q^2 um} \biggr\} d\tau.
\end{multline*} Hence,
$$
	V_{\chi, m, s} = \sum_{0 \les \nu \les \nu_1} \frac{c_{\nu} (\alpha)}{X^{2\nu}} h^{\delta (\nu + 1/2)} \frac{1}{q} \biggl( \frac{X^{1-\alpha}}{\alpha h} \biggr)^{\beta (\nu + 1/2)} \sum_{\sigma = -\infty}^{+\infty} \tau(\chi; \sigma) J(\sigma),
$$ where the meaning of $J(\sigma)$ is clear. We further simplify the last expression by setting
$$
	b_{\nu} (\alpha) = \frac{c_{\nu} (\alpha)}{\alpha^{\beta (\nu + 1/2)}}, 
$$ which gives
$$
	V_{\chi, m,s} = \frac{X}{q} \sum_{0 \les \nu \les \nu_1} b_{\nu} (\alpha) \biggl( \frac{h}{X^{\alpha}} \biggr)^{\nu + 1/2} \sum_{\sigma = -\infty}^{+\infty} \tau(\chi; \sigma) J(\sigma). 
$$ Let us denote
$$
	T_3 = \frac{(\alpha h q)^2 um}{4s X^{1-2\alpha}}, \qquad T_4 = 16 T_3 = \frac{4(\alpha h q)^2 um}{sX^{1-2\alpha}},
$$ and split the sum $V_{\chi, m, s}$ as follows:
\begin{multline*}
	V_{\chi, m, s} = \\ \frac{X}{q} \sum_{0 \les \nu \les \nu_1} b_{\nu} (\alpha) \biggl( \frac{h}{X^{\alpha}} \biggr)^{\nu + 1/2} \biggl( \sum_{T_3 < \sigma < T_4} + \sum_{|\sigma| \ges T_4} + \sum_{-T_4 < \sigma \les T_3} \biggr)
	\tau(\chi; \sigma) J(\sigma) =: \\ \frac{X}{q} \sum_{0 \les \nu \les \nu_1} b_{\nu} (\alpha) \biggl( \frac{h}{X^{\alpha}} \biggr)^{\nu + 1/2} \bigl( C_1 + C_2 + C_3 \bigr). 
\end{multline*} Similarly to above, we apply \hyperref[lemma2]{Lemma~2} to the integrals $J(\sigma)$ in $C_2$ and $C_3$ to estimate them from above and use \hyperref[lemma3]{Lemma~3} to compute $J(\sigma)$ in $C_1$. If $q$ is small enough and $T_4 < 1$, the whole sum $V_{\chi, m, s}$ is estimated by \hyperref[lemma2]{Lemma~2}. 

Next, we verify the conditions of \hyperref[lemma2]{Lemma~2} and~\hyperref[lemma3]{Lemma~3}. Put 
\begin{gather*}
	{\tilde w}(\tau) = f_2 \biggl( \frac{X^{1-\alpha} s\tau}{\alpha h qum} \biggr) \Psi_N\biggl( \frac{X^{1-\alpha} s\tau}{\alpha h qum} \biggr) \tau^{\beta (\nu + 1/2) - 1} \frac{d^{2\nu} G(t, \tau)}{dt^{2\nu}} \biggl|_{t=\tau^{\delta}}, \\
	{\tilde g}(\tau) = 
	\begin{cases}
	\displaystyle (1-\alpha) h X^{\alpha} \tau^{\gamma} - \frac{X^{1-\alpha} s\sigma \tau}{\alpha h q^2 um}, 
	&\text{if }
	\displaystyle \tau \in \bigl[ (1-\Delta)^{1/\delta}; (y+\Delta)^{1/\delta} \bigr],
	\\[2mm]
	\displaystyle 0
	&\text{if }
	\displaystyle \tau \les (1-2\Delta)^{1/\delta} \ \text{or } \tau \ges (y+2\Delta)^{1/\delta}
	\end{cases} 
\end{gather*} and define ${\tilde g} (\tau)$ on $(1-2\Delta)^{1-\alpha} \les \tau \les (1-\Delta)^{1-\alpha}$ and $(y+\Delta)^{1-\alpha} \les \tau \les (y+2\Delta)^{1-\alpha}$ appropriately.

Now we estimate ${\tilde w}^{(j)} (\tau)$ on $(1-\Delta)^{1-\alpha} \les \tau \les (y+\Delta)^{1-\alpha}$. We have
\begin{multline*}
	{\tilde w}^{(j)} (\tau) = \sum_{j_1+j_2+j_3+j_4 = j} \binom{j}{j_1, j_2, j_3, j_4} \frac{d^{j_1}}{d\tau^{j_1}} f_2 \biggl( \frac{X^{1-\alpha} s\tau}{\alpha h qum} \biggr) \frac{d^{j_2}}{d\tau^{j_2}} \Psi_N\biggl( \frac{X^{1-\alpha} s\tau}{\alpha h qum} \biggr) \cdot \\ \frac{d^{j_3}}{d\tau^{j_3}} \tau^{\beta(\nu + 1/2) - 1} \frac{d^{j_4}}{d\tau^{j_4}} \biggl( \frac{d^{2\nu} G(t, \tau)}{dt^{2\nu}} \biggl|_{t=\tau^{\delta}} \biggr).
\end{multline*} Next,
\begin{gather*}
	\frac{d^{j_1}}{d\tau^{j_1}} f_2 \biggl( \frac{X^{1-\alpha} s\tau}{\alpha h qum} \biggr) \ll \log X, \\
	\frac{d^{j_2}}{d\tau^{j_2}} \Psi_N\biggl( \frac{X^{1-\alpha} s\tau}{\alpha h qum} \biggr) \ll \biggl( \frac{X^{1-\alpha} s}{N\alpha hqum} \biggr)^{j_2} (\log X)^{A_0 j_2} \ll_{j_2} (\log X)^{A_0 j_2}, \\
	\frac{d^{j_3}}{d\tau^{j_3}} \bigl( \tau^{\beta (\nu + 1/2) - 1} \bigr) = \biggl( \beta\bigl(\nu + \frac{1}{2}\bigr) - 1 \biggr)_{j_3} \tau^{\beta (\nu + 1/2) - 1 - j_3} \ll_{j_3} 1.
\end{gather*} To estimate the last factor 
$$
	\frac{d^{j_4}}{d\tau^{j_4}} \biggl( \frac{d^{2\nu} G(t, \tau)}{dt^{2\nu}} \biggl|_{t=\tau^{\delta}} \biggr)
$$ we apply the formula of Faa di Bruno: for $r$ times differentiable functions $f(x)$, $\phi(x)$ one has
\begin{equation} \label{Faa}
	\frac{d^r \phi(f(x))}{dx^r} = \sum_{\substack{m_1 + 2m_2 + \ldots + rm_r = r \\ m_1, \ldots, m_r \ges 0}} \frac{r!}{m_1!\ldots m_r!} \phi^{(m_1 + \ldots + m_r)} (f(x)) \prod_{\kappa=1}^r \biggl( \frac{f^{(\kappa)}(x)}{\kappa!} \biggr)^{m_{\kappa}}.
\end{equation} First, let us compute $(d^{2\nu} / dt^{2\nu}) G(t, \tau)$. By the binomial theorem
$$
	\frac{d^{2\nu}}{dt^{2\nu}} G(t, \tau) = \sum_{l=0}^{2\nu} \binom{2\nu}{l} w^{(l)}(t) \frac{d^{2\nu - l}}{dt^{2\nu - l}} e(H(t, \tau)).
$$ Put in~\eqref{Faa} $f = H(t)$, $\phi = e(H)$. Then we get
\begin{multline} \label{first_Faa}
	\frac{d^{2\nu - l}}{dt^{2\nu - l}} e(H(t, \tau)) = \\ \sum_{\substack{m_1 + 2m_2 + \ldots + (2\nu-l)m_{2\nu-l} = 2\nu - l \\ m_1, \ldots, m_{2\nu-l} \ges 0}} \frac{(2\nu - l)!}{m_1! \ldots m_{2\nu-l}!} (2\pi i)^{m_1 + \ldots + m_{2\nu - l}} e(H(t, \tau)) \cdot \prod_{\kappa=1}^{2\nu - l} \biggl( \frac{H^{(\kappa)}(t, \tau)}{\kappa!} \biggr)^{m_{\kappa}}.
\end{multline} Note that the last expression is evaluated at $t = t_0 = \tau^{\delta}$. By definition,
\begin{equation} \label{def_H_again}
	H(t, \tau) = g(t, \tau) - g(t_0, \tau) - \frac{1}{2} g''(t_0, \tau) (t-t_0)^2,
\end{equation} and so $H(t_0, \tau) = H'(t_0, \tau) = H''(t_0, \tau) = 0$, $e(H(\tau^{\delta}, \tau)) = 1$. It means that the non-zero contribution to the right hand side of~\eqref{first_Faa} when $t = t_0$ would only come from the tuples of the form $(0, 0, m_3, m_4, \ldots)$.

Thus, we need to apply binomial theorem and Faa di Bruno's formula again to compute
$$
	\frac{d^{j_4}}{d\tau^{j_4}} \biggl( w^{(l)} (t) \biggr|_{t=\tau^{\delta}} \cdot \prod_{\kappa=3}^{2\nu - l} \biggl( \frac{1}{\kappa!} \frac{d^{\kappa} H(t, \tau)}{dt^{\kappa}} \biggr|_{t=\tau^{\delta}} \biggr)^{m_{\kappa}} \biggr).
$$ By the binomial theorem it is a linear combination of the expressions
\begin{equation} \label{two_factors}
	\frac{d^{j_5}}{d\tau^{j_5}} \biggl( w^{(l)} (t) \biggr|_{t=\tau^{\delta}} \biggr) \cdot \frac{d^{j_6}}{d\tau^{j_6}} \biggl( \prod_{\kappa=3}^{2\nu - l} \biggl( \frac{1}{\kappa!} \frac{d^{\kappa} H(t, \tau)}{dt^{\kappa}} \biggr|_{t=\tau^{\delta}} \biggr)^{m_{\kappa}} \biggr),
\end{equation} where $j_5 + j_6 = j_4$. First, we evaluate the factor
$$
	\frac{d^{j_5}}{d\tau^{j_5}} \biggl( w^{(l)} (t) \biggr|_{t=\tau^{\delta}} \biggr).
$$ Note that since $t_0 = t$, we have
$$
	\frac{d^r}{dt_0^r} w^{(l)} (t) \biggr|_{t=t_0} = \frac{d^{l+r}}{dt_0^{l+r}} w(t_0) = w^{(l+r)} (t_0).
$$ Then by Faa di Bruno's formula and~\eqref{bound_derivative_w}
\begin{multline*}
	\frac{d^{j_5}}{d\tau^{j_5}} \biggl( w^{(l)} (t) \biggr|_{t=\tau^{\delta}} \biggr) = \\ \sum_{\substack{l_1 + 2l_2 + \ldots + j_5 l_{j_5} = j_5 \\ l_1, \ldots, l_{j_5} \ges 0}} \frac{j_5 !}{l_1 ! \ldots l_{j_5}!} w^{(l+l_1+\ldots+l_{j_5})} (\tau^{\delta}) \prod_{\kappa=1}^{j_5} \biggl( \frac{d^{\kappa}}{d\tau^{\kappa} } t_0 (\tau) \biggr)^{\kappa} \ll_{j_5} \\  \sum_{\substack{l_1 + 2l_2 + \ldots + j_5 l_{j_5} = j_5 \\ l_1, \ldots, l_{j_5} \ges 0}} w^{(l+l_1+\ldots+l_{j_5})} (\tau^{\delta}) \ll w^{(l+j_5)} (\tau^{\delta}) \ll (\log X)^{C_0 (l+j_5) + 1}.
\end{multline*} Next, we evaluate the second factor in~\eqref{two_factors}. For $\kappa \ges 3$ we have
\begin{multline} \label{similarly_to}
	\frac{d^r}{d\tau^r} \biggl( \frac{1}{\kappa!} \frac{d^{\kappa} H(t, \tau)}{dt^{\kappa}} \biggr|_{t=\tau^{\delta}} \biggr) = \frac{1}{\kappa!} \frac{d^r}{d\tau^r} g^{(\kappa)} (\tau^{\delta}, \tau) = \frac{(\alpha)_{\kappa} h X^{\alpha}}{\kappa!} \frac{d^r}{d\tau^r} \tau^{\delta (\alpha - \kappa)} = \\ \frac{(\alpha)_{\kappa} h X^{\alpha}}{\kappa!} (\delta (\alpha - \kappa))_r \tau^{\delta (\alpha - \kappa) - r} \ll_{\kappa} hX^{\alpha}
\end{multline} and, thus,
$$
	\frac{d^{j_6}}{d\tau^{j_6}} \biggl( \prod_{\kappa=3}^{2\nu - l} \biggl( \frac{1}{\kappa!} \frac{d^{\kappa} H(t, \tau)}{dt^{\kappa}} \biggr|_{t=\tau^{\delta}} \biggr)^{m_{\kappa}} \biggr) \ll_{\kappa} (hX^{\alpha})^{m_3 + \ldots + m_{\kappa}} \ll (hX^{\alpha})^{(2\nu - l)/3},
$$ since $m_1 + 2m_2 + 3m_3 + \ldots + \kappa m_{\kappa} \ges m_3 + \ldots + m_{\kappa}$ and so the expression~\eqref{two_factors} can be bounded from above by 
$$
	(\log X)^{C_0 (l+j_5)+1} (hX^{\alpha})^{(2\nu - l)/3}.
$$ The last expression reaches its maximum at $l=0$, $j_5 = j_4$, $j_6 = 0$. We conclude
$$
	\frac{d^{j_4}}{d\tau^{j_4}} \biggl( \frac{d^{2\nu} G(t, \tau)}{dt^{2\nu}} \biggl|_{t=\tau^{\delta}} \biggr) \ll (\log X)^{C_0 j_4 + 1} (hX^{\alpha})^{2\nu / 3}.
$$ Finally, we find
\begin{multline} \label{w_1_derivative}
	{\tilde w}^{(j)}(\tau) \ll 
	(\log X)^2 (hX^{\alpha})^{2\nu  /3} \sum_{j_1 + \ldots + j_4 = j} \binom{j}{j_1, \ldots, j_4} (\log X)^{A_0 j_2 + C_0 j_4} \ll \\
	(\log X)^{2 + E_0 j} (hX^{\alpha})^{2\nu / 3},
\end{multline} where $E_0 = \max\{ A_0, C_0 \}$. So the inequality ${\tilde w}^{(j)} (\tau) \ll X_I Y_I^{-j}$ holds with $X_I = (\log X)^2 (hX^{\alpha})^{2\nu /3}$, $V_I = (\log X)^{-E_0}$. Next, we have
\begin{gather*}
	{\tilde g}'(\tau) = \gamma (1-\alpha) hX^{\alpha} \tau^{\gamma-1} - \frac{X^{1-\alpha} s\sigma}{\alpha h q^2 um}, \qquad
	{\tilde g}''(\tau) = \gamma (\gamma-1)(1-\alpha) hX^{\alpha} \tau^{\gamma-2}, \\
	{\tilde g}^{(j)} (\tau) = (\gamma)_j (1-\alpha) hX^{\alpha} \tau^{\gamma-j} \asymp hX^{\alpha},
\end{gather*} so one can take $Y_I = hX^{\alpha}, Q_I = 1$. Put $\mathbb{J} = \bigl[10^{-1}; 10\bigr]$, $\tilde V_I = |\mathbb{J}|$ and $Z_I = Q_I + X_I + Y_I + \tilde V_I + 1 \asymp (\log X)^2 (hX^{\alpha})^{2\nu / 3} + hX^{\alpha}$, which implies
$$
	Z_I  \asymp
	\begin{cases}
	\displaystyle hX^{\alpha},
	& \text{if } \displaystyle \nu = 0,1,
	\\[2mm]
	\displaystyle (\log X)^2 (hX^{\alpha})^{2\nu / 3},
	&\text{if } \displaystyle \nu \ges 2.
	\end{cases}
$$ We choose the constant $\delta_I > 0$ such that $Y_I > Z_I^{3\delta_I}$. If $\nu = 0,1$, we get $hX^{\alpha} > (hX^{\alpha})^{3\delta_I}$, which holds true for all $\delta_I < 1/3$. If $\nu \ges 2$, then
$$
	hX^{\alpha} \ges (\log X)^{6\delta_I} (hX^{\alpha})^{2\nu \delta_I},
$$ so one can take $\delta_I = 1 / (121 D_0)$. It is easy to check that
$$
	\frac{Q_I Z_I^{\delta_I / 2}}{\sqrt{Y_I}} \les V_I
$$ holds true for all $\nu \les \nu_1 = 60 D_0$.

Next, if $|\sigma| \ges T_4$, then
\begin{multline*}
	|{\tilde g}'(\tau)| = \biggl| \frac{X^{1-\alpha} s\sigma}{\alpha hq^2 um} - \gamma (1-\alpha) hX^{\alpha} \tau^{\gamma-1} \biggr| \ges \frac{X^{1-\alpha} s|\sigma|}{\alpha h q^2 um} \biggl( 1 - \frac{\alpha^2 h^2 q^2 um\tau^{\gamma-1}}{X^{1-2\alpha} sT_4} \biggr) \ges \\ \frac{2}{3} \frac{X^{1-\alpha} s|\sigma|}{\alpha hq^2 um}.
\end{multline*} If $-T_4 < \sigma \les 0$, then
$$
	{\tilde g}'(\tau) = \gamma(1-\alpha) hX^{\alpha} \tau^{\gamma-1} - \frac{X^{1-\alpha} s\sigma}{\alpha hq^2 um} = \alpha h X^{\alpha} \tau^{\gamma-1} + \frac{X^{1-\alpha} s|\sigma|}{\alpha hq^2 um} \ges \alpha hX^{\alpha} \tau^{\gamma-1} \ges \frac{1}{2} \alpha hX^{\alpha}.
$$ Finally, if $1 \les \sigma \les T_3$, then
$$
	{\tilde g}'(\tau) = \alpha hX^{\alpha} \tau^{\gamma-1} \biggl( 1- \frac{X^{1-2\alpha} s\sigma \tau^{1-\gamma}}{(\alpha hq)^2 um} \biggr) \ges \alpha hX^{\alpha} \tau^{\gamma-1} \biggl( 1-\frac{3}{5} \biggr) \ges \frac{1}{6} \alpha h X^{\alpha}.
$$ So one can choose
$$
	R_I = 
	\begin{cases}
	\displaystyle \frac{2}{3} \frac{X^{1-\alpha} s|\sigma|}{\alpha h q^2 um} 
	&\text{if } \displaystyle |\sigma| \ges T_4,
	\\[3mm]
	\displaystyle \frac{1}{6} \alpha h X^{\alpha}
	&\text{if } \displaystyle -T_4 < \sigma \les T_3.
	\end{cases}
$$ Again, setting
$$
	\Delta_1 = \frac{Q_I R_I}{\sqrt{Y_I}}, \qquad \Delta_2 = R_I V_I,
$$ we show that $\Delta_1, \Delta_2 > 1$. Indeed, in the case $|\sigma| \ges T_4$, we have
\begin{gather*}
	\Delta_1 = \frac{2}{3} \frac{X^{1-\alpha} s|\sigma|}{\alpha h q^2 um} \frac{1}{\sqrt{hX^{\alpha}}} \ges \frac{2}{3\sqrt{h}} \frac{X^{1-3\alpha/2} s}{\alpha h q^2 um} T_4 = \frac{8}{3} \alpha \sqrt{h} X^{\alpha/2} > 1, \\
	\Delta_2 = \frac{2}{3} \frac{X^{1-\alpha} s|\sigma|}{\alpha h q^2 um} (\log X)^{-E_0} \ges \frac{2}{3} \frac{X^{1-\alpha} s}{\alpha h q^2 um} T_4 (\log X)^{-E_0} = \frac{8}{3} \alpha \bigl(hX^{\alpha}\bigr)^{2/3} (\log X)^{-E_0} > 1.
\end{gather*} If $-T_4 < \sigma \les T_3$, then
\begin{gather*}
	\Delta_1 = \frac{1}{6} \alpha h X^{\alpha} \frac{1}{\sqrt{hX^{\alpha}}} = \frac{\alpha}{6} \sqrt{h} X^{\alpha/2} > 1, \\
	\Delta_2 = \frac{1}{6} \alpha hX^{\alpha} (\log X)^{-E_0} > 1.
\end{gather*} Thus, applying \hyperref[lemma2]{Lemma~2} with a large enough $F_0 = F_0 (\alpha) > 1$, $\mathbb{J} = \bigl[ (1-\Delta)^{1/\delta}; (y + \Delta)^{1/\delta} \bigr]$, for $|\sigma| \ges T_4$ we find
\begin{multline*}
	J(\sigma) \ll |\mathbb{J}| X_I (\Delta_1^{-F_0} + \Delta_2^{-F_0}) \ll \\ (\log X)^2 (hX^{\alpha})^{2\nu / 3} \biggl\{ \biggl( \frac{\alpha h q^2 um \sqrt{hX^{\alpha}}}{2X^{1-\alpha} s|\sigma|} \biggr)^{F_0} + \biggl( \frac{3\alpha h q^2 um (\log X)^{E_0}}{2X^{1-\alpha} s|\sigma|} \biggr)^{F_0} \biggr\} \ll \\ (\log X)^2 (hX^{\alpha})^{2\nu/3} \biggl( \frac{3\alpha}{2} h^{3/2} q^2 um \frac{1}{X^{1-3\alpha/2} s |\sigma|} \biggr)^{F_0},
\end{multline*} and for $-T_4 < \sigma \les T_3$ we get
\begin{multline*}
	J(\sigma) \ll (\log X)^2 (hX^{\alpha})^{2\nu/3} \biggl\{ \biggl( \frac{6}{\alpha \sqrt{h} X^{\alpha/2}} \biggr)^{F_0} + \biggl( \frac{6 (\log X)^{E_0}}{\alpha h X^{\alpha}} \biggr)^{F_0} \biggr\} \\ \ll (\log X)^2 (hX^{\alpha})^{2\nu/3} \biggl( \frac{6}{\alpha \sqrt{h}} X^{-\alpha/2} \biggr)^{F_0}.
\end{multline*} Thus, the contribution from $C_2$, $C_3$ to $V_{\chi, m, s}$ can be made small enough with the appropriate choice of $F_0$. We get the formula
\begin{gather} \label{C_2_C_3_contribution}
	V_{\chi,m,s} = \frac{X}{q} \sum_{0 \les \nu \les \nu_1} b_{\nu} (\alpha) \biggl( \frac{h}{X^{\alpha}} \biggr)^{\nu + 1/2} C_1 + O\bigl(X^{- \alpha F_0 / 10}\bigr), \\
	C_1 = \sum_{T_3 < \sigma < T_4} \tau(\chi; \sigma) J(\sigma). \nonumber 
\end{gather} We are now ready to compute $C_1$ using \hyperref[lemma3]{Lemma~3}. Let $T_3 < \sigma < T_4$, ${\tilde g}'(\tau_0) = 0$. Then $\tau_0 \in \mathbb{J} = \bigl[10^{-1}; 10\bigr]$. We find
$$
	\tau_0 = \frac{1}{X^{1-\alpha}} \biggl\{ \frac{(\alpha h q)^2 um}{s\sigma} \biggr\}^{\xi},
$$ where $\xi = 1/(1-\gamma) = (1-\alpha) / (1-2\alpha)$; 
$$
	{\tilde g} (\tau_0) = (1-\alpha) hX^{\alpha} \tau_0^{\gamma} - \frac{X^{1-\alpha} s\sigma \tau_0}{\alpha hq^2 um} = (1-2\alpha) h \biggl\{ \frac{(\alpha h q)^2 um}{s\sigma} \biggr\}^{\eta},
$$ where $\eta = \alpha / (1-2\alpha)$;
$$
	{\tilde g}''(\tau_0) = -\frac{\alpha (1-2\alpha)}{1-\alpha} hX^{2(1-\alpha)} \biggl( \frac{s\sigma}{(\alpha h q)^2 um} \biggr)^{\omega},
$$ where $\omega = \xi (2-\gamma) = (2-3\alpha) / (1-2\alpha)$. Finally, take
$$
	{\tilde G} (\tau) = {\tilde w} (\tau) e^{2\pi i {\tilde H} (\tau)}, \qquad {\tilde H} (\tau) = {\tilde g} (\tau) - {\tilde g} (\tau_0) - \frac{{\tilde g}''(\tau_0)}{2} (\tau - \tau_0)^2.
$$ Then \hyperref[lemma3]{Lemma~3} implies
\begin{multline*}
	J(\sigma) = e\biggl( {\tilde g} (\tau_0) - \frac{1}{8} \biggr) \sum_{0 \les \mu \les \mu_1} \frac{\sqrt{2\pi}}{\mu!} \frac{(2i)^{-\mu}}{|{\tilde g}''(\tau_0)|^{\mu + 1/2}} \frac{d^{2\mu}  {\tilde G} (\tau)}{d\tau^{2\mu}} \biggr|_{\tau=\tau_0} + O\bigl(X^{-\alpha F_0}\bigr) = \\ e \biggl( (1-2\alpha) h \biggl\{ \frac{(\alpha h q)^2 um}{s\sigma} \biggr\}^{\eta} - \frac{1}{8} \biggr) \sum_{0 \les \mu \les \mu_1} \frac{\sqrt{2\pi}}{\mu!} \frac{(2i)^{-\mu}}{\bigl( \alpha (1-2\alpha) / (1-\alpha) \bigr)^{\mu + 1/2}} \cdot \\  \frac{1}{(hX^{2(1-\alpha)})^{\mu+1/2}} \biggl( \frac{(\alpha h q)^2 um}{s\sigma} \biggr)^{\omega(\mu+1/2)} \frac{d^{2\mu}  {\tilde G} (\tau)}{d\tau^{2\mu}} \biggr|_{\tau=\tau_0} + O\bigl(X^{-\alpha F_0}\bigr)
\end{multline*} with $\mu_1 = 3F_0 / \delta_I = 363D_0 F_0$. Setting
$$
	c_{\mu} (\alpha) = e\biggl( \frac{1}{8} \biggr) \frac{\sqrt{2\pi}}{\mu!} (2i)^{-\mu} \biggl( \frac{1-\alpha}{\alpha (1-2\alpha)} \biggr)^{\mu + 1/2} \alpha^{2\omega(\mu +1/2)},
$$ we get
\begin{multline*}
	J(\sigma) = e \biggl( (1-2\alpha) h \biggl\{ \frac{(\alpha h q)^2 um}{s\sigma} \biggr\}^{\eta} \biggr) \sum_{0 \les \mu \les \mu_1} \frac{c_{\mu} (\alpha)}{\bigl(hX^{2(1-\alpha)}\bigr)^{\mu + 1/2}} \biggl( \frac{(hq)^2 um}{s\sigma} \biggr)^{\omega(\mu + 1/2)} \cdot \\ \frac{d^{2\mu}  {\tilde G} (\tau)}{d\tau^{2\mu}} \biggr|_{\tau=\tau_0} + O\bigl(X^{-\alpha F_0}\bigr).
\end{multline*} Again, it is not hard to see that the $O$-term contributes at most $O\bigl(  X^{-\alpha F_0 / 10}\bigr)$ to the sum $V_{\chi, m, s}$. This contribution can be made arbitrarily small. Hence, from~\eqref{C_2_C_3_contribution} we have
\begin{multline*}
	V_{\chi, m,s} = \\ \frac{X}{q} \sum_{T_3 < \sigma \les T_4} \tau(\chi; \sigma) \sum_{0\les \nu \les \nu_1} \sum_{0\les \mu \les \mu_1} b_{\nu} (\alpha) \biggl( \frac{h}{X^{\alpha}} \biggr)^{\nu + 1/2} \frac{c_{\mu} (\alpha)}{\bigl(hX^{2(1-\alpha)}\bigr)^{\mu + 1/2}} \biggl( \frac{(hq)^2 um}{s\sigma} \biggr)^{\omega(\mu + 1/2)} \cdot \\ \frac{d^{2\mu}  {\tilde G} (\tau)}{d\tau^{2\mu}} \biggr|_{\tau=\tau_0} e\biggl( (1-2\alpha) h \biggl\{ \frac{(\alpha h q)^2 um}{s\sigma} \biggr\}^{\eta} \biggr) + O\bigl(X^{- \alpha F_0 / 10} \bigr).
\end{multline*} Substituting this expression into~\eqref{V_chi_intro} and changing the order of summation, we get
\begin{multline} \label{W_M_N_K_new}
	W(M,N,K) = \\ \frac{X}{q^2} \sum_{m=1}^{+\infty} \frac{f_1 (m)}{m} \Psi_M(m) \sum_{\substack{T_1 < s < T_2 \\ T_3 < \sigma < T_4}} \sum_{\substack{0 \les \nu \les \nu_1 \\ 0 \les \mu \les \mu_1}} b_{\nu} (\alpha) c_{\mu} (\alpha) \biggl( \frac{h}{X^{\alpha}} \biggr)^{\nu + 1/2} \frac{1}{\bigl(hX^{2 (1-\alpha)}\bigr)^{\mu + 1/2}} \cdot \\ \biggl( \frac{(hq)^2 um}{s\sigma} \biggr)^{\omega (\mu + 1/2)} \frac{d^{2\mu}  {\tilde G} (\tau)}{d\tau^{2\mu}} \biggr|_{\tau=\tau_0} e\biggl( (1-2\alpha) h \biggl\{ \frac{(\alpha h q)^2 um}{s\sigma} \biggr\}^{\eta} \biggr) \cdot \\ \frac{1}{\varphi(q)} \sum_{\chi \ \text{mod} \ q} \chi(mua^{\ast}) \tau(\chi; s) \tau(\chi; \sigma) + O\bigl(X^{-\alpha F_0/20}\bigr) + O\bigl(X^{-\alpha D_0 /2}\bigr).
\end{multline}

\subsection*{Bounding Kloosterman sum}

We rewrite the inner sum in~\eqref{W_M_N_K_new} as follows:
\begin{multline} \label{equation_smth}
\frac{1}{\varphi(q)} \sum_{\chi \ \text{mod} \ q} \chi(mua^{\ast}) \tau(\chi; s) \tau(\chi; \sigma) = \sum_{l,r=1}^q e\biggl( \frac{ls + r\sigma}{q} \biggr) \frac{1}{\varphi(q)} \sum_{\chi \ \text{mod} \ q} \chi(lrmua^{\ast}).
\end{multline} By orthogonality of characters the sum in the right hand side of~\eqref{equation_smth} transforms into Kloosterman sum
$$
\sum_{\substack{l=1 \\ (l,q)=1}}^q e\biggl( \frac{sl + \sigma a(mu)^{\ast} l^{\ast}}{q} \biggr) = S_q \bigl(s, \sigma a(mu)^{\ast}\bigr).
$$ Thus,
\begin{multline*}
	W(M,N,K) = \frac{X}{uq^2} \sum_{\substack{m=1 \\ (m,q)=1}}^{\infty} \frac{f_1 (m)}{m} \Psi_M(m) \sum_{\substack{T_1 < s < T_2 \\ T_3 < \sigma < T_4}} \sum_{\substack{0 \les \nu \les \nu_1 \\ 0 \les \mu \les \mu_1}} b_{\nu} (\alpha) c_{\mu} (\alpha) \biggl( \frac{h}{X^{\alpha}} \biggr)^{\nu + 1/2} \cdot \\ \frac{1}{\bigl( hX^{2(1-\alpha)} \bigr)^{\mu + 1/2}} \biggl( \frac{(hq)^2 um}{s\sigma} \biggr)^{\omega (\mu+1/2)} \frac{d^{2\mu}  {\tilde G} (\tau)}{d\tau^{2\mu}} \biggr|_{\tau=\tau_0} \cdot \\ e\biggl( (1-2\alpha) h \biggl\{ \frac{(\alpha h q)^2 um}{s\sigma} \biggr\}^{\eta} \biggr) S_q \bigl(s, \sigma a(mu)^{\ast}\bigr) + O\bigl(X^{-\alpha F_0/20}\bigr) + O\bigl(X^{- \alpha D_0/2}\bigr).
\end{multline*} We can now estimate the multiple sum over $m,s,\sigma,v$ and $\mu$. Since ${\tilde H} (\tau_0) = {\tilde H}'(\tau_0) = {\tilde H}''(\tau_0) = 0$, similarly to~\eqref{similarly_to}, we get the upper bound
$$
	\frac{d^r}{d\tau^r} \bigl( e\bigl( {\tilde H} (\tau) \bigr) \bigr) \biggr|_{\tau=\tau_0} \ll (hX^{\alpha})^{r/3}.
$$ Together with~\eqref{w_1_derivative} this implies
$$
{\tilde G}^{(2\mu)} (\tau_0) \ll (hX^{\alpha})^{(2/3) (\nu + \mu)} (\log X)^2.
$$ Next, we apply Weil's bound:
$$
\bigl|S_q (s, \sigma a (mu)^{\ast})\bigr| \les \tau(q) \sqrt{q} (s,\sigma, q)^{1/2}.
$$ Changing the order of summation, we get the inequality
\begin{multline*}
W(M,N,K) \ll \frac{X \tau(q)}{uq\sqrt{q}} \frac{(\log X)^3}{M} \sum_{\substack{0 \les \nu \les \nu_1 \\ 0 \les \mu \les \mu_1}} \biggl( \frac{h}{X^{\alpha}} \biggr)^{\nu + 1/2}  \frac{(hX^{\alpha})^{(2/3) (\nu + \mu)}}{\bigl(hX^{2(1-\alpha)}\bigr)^{\mu + 1/2}} (hq)^{\omega (2\mu+1)} \cdot \\ \sum_{M/\Theta \les m \les M\Theta} (mu)^{\omega (\mu + 1/2)} \sum_{\substack{T_1 < s < T_2 \\ T_3 < \sigma < T_4}} \frac{(s, \sigma, q)^{1/2}}{(s\sigma)^{\omega (\mu + 1/2)}} +  O\bigl(X^{-\alpha F_0/20} + X^{- \alpha D_0/2}\bigr).
\end{multline*} The sums over $s$ and $\sigma$ could be bounded as
$$
\sum_{T_1 < s < T_2} \sum_{T_3 < \sigma < T_4} \frac{(s, \sigma, q)^{1/2}}{(s\sigma)^{\omega (\mu + 1/2)}} \ll_{\alpha} \biggl( \frac{X^{1-2\alpha}}{(hq)^2 um} \biggr)^{\omega (\mu + 1/2)} \sum_{T_1 < s < T_2} \sum_{T_3 < \sigma < T_4} (s,\sigma)^{1/2}.
$$ The last expression does not exceed
\begin{multline*}
\biggl( \frac{X^{1-2\alpha}}{(hq)^2 um} \biggr)^{\omega (\mu + 1/2)} \sum_{1 \les d \les \min(T_2, 16 T_4)} \sum_{\substack{T_1 < s < T_2 \\ s \equiv 0\pmod{d}}} \sum_{\substack{T_3 / 16 < \sigma < 16 T_4 \\ \sigma \equiv 0\pmod{d}}} \sqrt{d} \ll \\ \biggl( \frac{X^{1-2\alpha}}{(hq)^2 um} \biggr)^{\omega(\mu + 1/2)} \sum_{1 \les d \les \min(T_2, 16 T_4)} \sqrt{d} \frac{T_2}{d} \frac{T_4}{d} \ll \\ T_2 T_6 \biggl( \frac{X^{1-2\alpha}}{(hq)^2 um} \biggr)^{\omega (\mu + 1/2)} \ll \biggl( \frac{X^{1-2\alpha}}{(hq)^2 um} \biggr)^{\omega (\mu + 1/2) - 1}.
\end{multline*} Next, the summation over $m$ gives
\begin{multline*}
\sum_{M/\Theta \les m \les M \Theta} (mu)^{\omega (\mu + 1/2)} \cdot \biggl( \frac{X^{1-2\alpha}}{(hq)^2 mu} \biggr)^{\omega (\mu +1/2) - 1} \ll \\ M^2 U \biggl( \frac{X^{1-2\alpha}}{(hq)^2} \biggr)^{\omega (\mu + 1/2) - 1} \frac{1}{(\log X)^{A_0}},
\end{multline*} whence, if $D_0$ and $F_0$ are sufficiently large, 
\begin{multline*}
W(M, N, K) \ll \frac{X \tau(q)}{Uq\sqrt{q}} \frac{(\log X)^3}{M} \frac{M^2 U}{(\log X)^{A_0}} \sum_{\substack{0 \les \nu \les \nu_1 \\ 0 \les \mu \les \mu_1}} \biggl( \frac{h}{X^{\alpha}} \biggr)^{\nu + 1/2} \frac{(hX^{\alpha})^{(2/3)(\nu + \mu)}}{\bigl(hX^{2(1-\alpha)}\bigr)^{\mu + 1/2}} \cdot \\ (hq)^{\omega (2\mu + 1)} \biggl( \frac{X^{1-2\alpha}}{(hq)^2} \biggr)^{\omega (\mu +1/2) - 1} \ll \frac{X \tau(q)}{q \sqrt{q}} (\log X)^{3-A_0} M \sum_{\substack{0 \les \nu \les \nu_1 \\ 0 \les \mu \les \mu_1}} X^{\kappa_1} h^{\kappa_2} q^2,
\end{multline*} where
$$
\kappa_1 = -\frac{\alpha \nu}{3} - \frac{\alpha \mu}{3} + \alpha - 1, \qquad
\kappa_2 = 2 + \frac{5\nu}{3} - \frac{\mu}{3}.
$$ Clearly the main contribution comes from the term $\nu = \mu = 0$. We get
$$
W(M,N,K) \ll X^{\alpha} \sqrt{q} \tau(q) (\log X)^{3-A_0} M h^2.
$$ Summing $W(M,N,K)$ over all admissible $U,u,M,N,K$, and using Mardzhanishvili's inequality
$$
\sum_{u \les 2U} \tau_7 (u) \ll U (\log U)^6,
$$ finally, we find
\begin{multline} \label{contribution_type_III}
	W_{III} \ll X^{\alpha} \sqrt{q} \tau(q) (\log X)^{2C + 3 - A_0} \mathop{{\sum}^{'}}_{U \in \mathbf{G}} \sum_{U\Theta^{-7} \les u \les U\Theta^7} \bigl|F(u, U)\bigr| \cdot \\ \mathop{{\sum}^{'}}_{\substack{M_1 \les M \les M_2 \\ M \in \mathbf{G}}} M \mathop{{\sum}^{'}}_{\substack{N_1 \les N \les N_2 \\ N \in \mathbf{G}}} \mathop{{\sum}^{'}}_{\substack{K_1 \les K \les K_2 \\ K \in \mathbf{G}}} 1 \ll \\ X^{\alpha} \sqrt{q} \tau(q) (\log X)^{2C + 3 - A_0} \mathop{{\sum}^{'}}_{U \in \mathbf{G}} \sum_{u \les 2U} \tau_7 (u) (\log X)^{6(A_0 + 1) + 1} \biggl( \frac{X}{U} \biggr)^{1/3} (\log X)^{3(A_0 + 1)} \ll \\ X^{1/3 + \alpha} \sqrt{q} \tau(q) (\log X)^{2C + 3 - A_0 + 9(A_0 + 1) + 1} \mathop{{\sum}^{'}}_{U \in \mathbf{G}} U^{2/3} (\log U)^6 \ll \\ X^{1/3 + \alpha} \sqrt{q} \tau(q) X^{(2/3)(1/10 - 3\varepsilon_1 / 2)} (\log X)^{8A_0 + 2C + 19} \ll \\ X^{2/5 + \alpha - \varepsilon_1} \sqrt{q} \tau(q) (\log X)^{L_0},
\end{multline} where $L_0 = 8A_0 + 2C + 19$.

\subsection*{Final bound}

Combining together type~I~\eqref{contribution_type_I}, type~II~\eqref{contribution_type_II} and type~III~\eqref{contribution_type_III} estimates we get
\begin{multline*}
W \ll X^{2/5 + \alpha / 2 - \varepsilon_1 + 2\delta_1} + \frac{X^{1 - \alpha/2 + 2\delta_1}}{q} +
\frac{1}{\sqrt{q}} X^{4/5 + \alpha/4 + \varepsilon_1 / 2 + 3\delta_2} + \\ \frac{1}{q} X^{1 - \alpha/4 + 3\delta_2} +
X^{2/5 + \alpha - \varepsilon_1} \sqrt{q} \tau(q) (\log X)^{L_0}.
\end{multline*} Further, the right hand side of the last inequality does not exceed
\begin{multline*}
\frac{X}{q} \biggl( q X^{-3/5 + \alpha/2 - \varepsilon_1 + 2\delta_1} + X^{-\alpha/2 + 2\delta_1} + \sqrt{q} X^{-1/5 + \alpha/4 + \varepsilon_1 / 2 +3\delta_2} + \\ X^{-\alpha/4 + 3\delta_2} + q^{3/2} X^{-3/5 + \alpha - \varepsilon_1 + \delta_4} \biggr)
\end{multline*} with an arbitrarily small $\delta_4 > 0$. Clearly
$$
\max \bigl( X^{-\alpha/2 + 2\delta_1}, X^{-\alpha/4 + 3\delta_2} \bigr) \ll (\log X)^{-A}
$$ if $\delta_1$ and $\delta_2$ are small enough. Next,
$$
q X^{-3/5 + \alpha/2 - \varepsilon_1 + 2\delta_1} \ll X^{-1/5 - \alpha/10 + 2\delta_1} \ll (\log X)^{-A}.
$$ Then
$$
W \ll \frac{X}{q} \biggl( (\log X)^{-A} + \max \bigl( \sqrt{q} X^{-1/5 + \alpha/4 + \varepsilon_1 / 2 +3\delta_2}, q^{3/2} X^{-3/5 + \alpha - \varepsilon_1 + \delta_4} \bigr) \biggr).
$$ Thus, $W \ll (X/q) (\log X)^{-A}$ if
$$
q \les \min \biggl( (\log X)^{-2A} X^{2/5 - \alpha/2 - \varepsilon_1 - 6\delta_2}, (\log X)^{-(2/3)A} X^{2/5 - (2/3)\alpha + (2/3) \varepsilon_1 - (2/3) \delta_4} \biggr).
$$ The maximum of this bound is reached at $\varepsilon_1 = \alpha/10$. Thus, $q \les X^{2/5 - (3/5)\alpha - \varepsilon}$ with any $\varepsilon < \min(6\delta_2, (2/3)\delta_4)$. Finally, the desired bound~\eqref{to_prove} follows from partial summation.

\begin{remark}
	One can obtain a slightly better level of distribution, $q \les X^{2/5 - \alpha/2 - \varepsilon}$, in~\hyperref[thm1]{Theorem~1} by iterating the Poisson summation for the third time (on the sum over~$m$) and applying the bound for 2-dimensional Kloosterman sum~\cite{24}.
\end{remark}

\section*{ACKNOWLEDGMENTS}
I would like to thank M.A.~Korolev and M.~Radziwill for many helpful discussions and advice. I also thank the anonymous referee for the careful reading of the paper and for pointing out a mistake in the earlier version.


\begin{thebibliography}{99}
	
	
	
	
	
	\bibitem{1}
	R.~Baker, G.~Kolesnik
	``On the distribution of $p^{\alpha}$ modulo one'',
	J. Reine Angew. Math.
	\textbf{1985}
	(356),
	174--193
	(1985).
	
	
	\bibitem{2}
	A.~Balog
	``On the fractional parts of $p^{\theta}$'',
	Arch. Math. (Basel)
	\textbf{40}
	(5),
	434--440
	(1983).
	
		
	\bibitem{3}
	V.~Blomer, R.~Khan, M.~Young,
	``Distribution of mass of holomorphic cusp forms'',
	Duke Math. J.
	\textbf{162} (14),
	2609--2644 (2013).

		
	\bibitem{4}
	X.~Cao, W.~Zhai
	``On the distribution of $p^{\alpha}$ modulo one'',
	J. Theor. Nombres Bordeaux.
	\textbf{11}
	(2),
	407--423
	(1999).
	
	
	\bibitem{5}
	M.~E.~Changa
	``Primes in special intervals and additive problems with such numbers'',
	Mathematical Notes
	\textbf{73}
	(3),
	423--436
	(2003).
	
	
	\bibitem{6}
	T.~Estermann,
	``On Kloosterman's sum'',
	Mathematika
	\textbf{8},
	83--86 (1961).
	
	
	\bibitem{7}
	E.~Fourvy,
	``Autour du theoreme de Bombieri-Vinogradov'',
	Acta Math.
	\textbf{152} (3-4),
	219--244 (1984).
	
	
		
	\bibitem{8}
	E.~Fourvy, H.~Iwaniec,
	``Primes in arithmetic progressions'',
	Acta Arith.
	\textbf{42} (2),
	197--218 (1983).
	
		
	\bibitem{9}
	E.~P.~Golubeva, O.~M.~Fomenko
	``On the distribution of the sequence $\{ bp^{3/2} \}$ mod 1'',
	J. Soviet Math.,
	\textbf{17}
	(5),
	2102--2107
	(1981).
	
	\bibitem{10}
	S.~A.~Gritsenko
	``On the problem of I.M.Vinogradov (in russian)'',
	Mat. Zametki
	\textbf{39}
	(5),
	625--640
	(1986).
	
	
	\bibitem{11}
	S.~A.~Gritsenko, N.~A.~Zinchenko
	``On a trigonometric sum over primes (in russian)'',
	Nauch. Ved. BelSU. Ser.: Mat. Fiz.
	\textbf{30}
	(5(148)),
	48--52
	(2013). 
	
	
	\bibitem{12}
	G.~Harman
	``On the distribution of $\sqrt{p}$ modulo one'',
	Mathematika
	\textbf{30}
	(1),
	104--116
	(1983).
	
	
	\bibitem{13}
	G.~Harman, P.~Lewis
	``Gaussian primes in narrow sectors'',
	Mathematika
	\textbf{48} (1-2),
	119--135
	(2003).
	
	
		
	\bibitem{14}
	D.~R.~Heath-Brown,
	``Prime numbers in short intervals and a generalized Vaughan identity'',
	Canadian J. Math.
	\textbf{34} (6),
	1365--1377 (1982).
	
	
		
	\bibitem{15}
	H.~Iwaniec, E.~Kowalski,
	\textit{Analytic Number Theory},
	American Mathematical Soc.,
	\textbf{53} (2004).
	
	
		
	\bibitem{16}
	A.~A.~Karatsuba
	\textit{Basic analytic number theory}, 
	Berlin: Springer-Verlag, 1993,
	2nd ed.
	
	
	
	\bibitem{17}
	R.~M.~Kaufman
	``The distribution of $\{\sqrt{p} \}$ (in russian)'',
	Mat. Zametki
	\textbf{26}
	(4),
	497--504
	(1979).
	
		
	\bibitem{18}
	D.~Leitmann
	``On the uniform distribution of some sequences'',
	J. London Math. Soc. (2)
	\textbf{14}
	(3),
	430--432
	(1976).
	
	
	\bibitem{19}
	Yu.~V.~Linnik
	``On a theorem in prime number theory (in russian)'',
	Dokl. Akad. Nauk SSSR
	\textbf{47}
	(1),
	7--8
	(1945).
	
	
		
	\bibitem{20}
	K.~K.~Mardzhanishvili
	``An estimate of one arithmetic sum'',
	Dokl. Akad. nauk SSSR
	\textbf{22}
	(7),
	391--393
	(1939). 
	
	
	\bibitem{21}
	D~.H.~J.~Polymath,
	``New equidistribution estimates of Zhang type'',
	Algebra Number Theory
	\textbf{8} (9),
	2067--2199 (2014).
	
	
		
	\bibitem{22}
	X.~M.~Ren
	``Vinogradov's exponential sum over primes'',
	Acta Arith.
	\textbf{124}
	(3),
	269--285
	(2006).
	
	
	\bibitem{23}
	A.~V.~Shubin,
	``Fractional parts of non-integer powers of primes'',
	Math. Notes 
	\textbf{108} (3), 2020;
	\href{arxiv: 2010.15216}{arxiv:~2010.15216} 
	
	
	\bibitem{24}
	R.~A.~Smith,
	``On n-dimensional Kloosterman sums'',
	J.~Number Theory
	\textbf{11} (3),
	324--343 (1979).
	
	
		
	\bibitem{25}
	D.~I.~Tolev
	``On a theorem of Bombieri-Vinogradov type for prime numbers from a thin set'',
	Acta Arith.
	\textbf{81}
	(1),
	57--68
	(1997). 
	
	
	
	\bibitem{26}
	I.~M.~Vinogradov
	``A general property of prime number distribution (in russian)'', 
	Mat. Sb.
	\textbf{7(49)}
	(2),
	365--372
	(1940).


	\bibitem{27}
	I.~M.~Vinogradov
	``On a trigonometric sum over primes (in russian)'',
	Izv. Akad. Nauk SSSR. Ser. Mat. 
	\textbf{23}
	(2),
	157--164
	(1959).
	
	
		
	\bibitem{28}
	I.~M.~Vinogradov,
	\textit{Special Versions of the Method of Trigonometric Sums} (in russian),
	Nauka, Moscow, 1976.
	
	
	
\end{thebibliography}
\end{document}